\documentclass[a4paper]{amsart}
\usepackage{amssymb} 
\usepackage[T1]{fontenc}
\usepackage[latin1]{inputenc}
\usepackage{amsfonts}
\usepackage{amsxtra}
\usepackage{ae}
\usepackage[all]{xy}
\usepackage{enumerate}

\include{diagram}

\newcommand*{\ket}{\rangle}
\newcommand*{\bra}{\langle}
\newcommand*{\ad}{\mathsf{ad}}

\newcommand*{\A}{\mathcal{A}}

\newcommand*{\E}{\mathcal{E}}
\newcommand*{\F}{\mathcal{F}}

\renewcommand*{\H}{\mathcal{H}}

\renewcommand*{\S}{\mathcal{S}}
\newcommand*{\CC}{\mathcal{CC}}
\newcommand*{\CI}{\mathcal{CI}}
\renewcommand*{\P}{\mathcal{P}}
\newcommand*{\Q}{\mathcal{Q}}

\newcommand*{\T}{\mathcal{T}}

\newcommand*{\twisted}{\boxtimes}

\newcommand*{\DD}{\mathsf{D}}

\newcommand*{\Rep}{\mathsf{Rep}}

\renewcommand*{\max}{\mathsf{f}}
\newcommand*{\red}{\mathsf{r}}

\newcommand*{\cop}{\mathsf{cop}}
\renewcommand*{\top}{\mathsf{top}}
\newcommand*{\Alg}{$-$\mathsf{Alg}}
\newcommand*{\hit}{\rightharpoonup}
\newcommand*{\hitby}{\leftharpoonup}

\newcommand*{\KH}{\mathbb{K}}
\newcommand*{\LH}{\mathbb{L}}

\DeclareMathOperator{\sgn}{sgn}
\DeclareMathOperator{\Mor}{Mor}

\DeclareMathOperator{\tr}{tr}

\DeclareMathOperator{\res}{res}
\DeclareMathOperator{\ind}{ind}

\DeclareMathOperator{\id}{id}

\numberwithin{equation}{section}
\theoremstyle{change}
\newtheorem{theorem}{Theorem}[section]
\newtheorem{prop}[theorem]{Proposition}
\newtheorem{lemma}[theorem]{Lemma}
\newtheorem{cor}[theorem]{Corollary}
\newtheorem{definition}[theorem]{Definition}

\begin{document}

\title[Baum-Connes conjecture]{The Baum-Connes conjecture for free orthogonal quantum groups}

\author{Christian Voigt}

\address{Christian Voigt \\
         Mathematisches Institut\\
         Westf\"alische Wilhelms-Universit\"at M\"unster \\
         Einsteinstra\ss e 62 \\
         48149 M\"unster\\
         Germany 
}
\email{cvoigt@math.uni-muenster.de}


\maketitle

\begin{abstract}
We prove an analogue of the Baum-Connes conjecture for free orthogonal quantum groups. 
More precisely, we show that these quantum groups have a 
$ \gamma $-element and that $ \gamma = 1 $. It follows that free orthogonal quantum groups 
are $ K $-amenable. We compute explicitly their $ K $-theory and deduce in the unimodular case that the 
corresponding reduced $ C^* $-algebras do not contain nontrivial idempotents. \\
Our approach is based on the reformulation of the Baum-Connes conjecture by Meyer and Nest using the language of
triangulated categories. An important ingredient 
is the theory of monoidal equivalence of compact quantum groups developed by Bichon, De Rijdt and Vaes. 
This allows us to study the problem in terms of the quantum group $ SU_q(2) $. 
The crucial part of the argument is a detailed analysis of the equivariant Kasparov theory of 
the standard Podle\'s sphere. 
\end{abstract}

\section{Introduction}

Let $ G $ be a second countable locally compact group and let $ A $ be a separable $ G $-$ C^* $-algebra. 
The Baum-Connes conjecture with coefficients in $ A $ asserts that the 
assembly map 
$$ 
\mu_A: K^{\top}_*(G;A) \rightarrow K_*(G \ltimes_\red A)
$$
is an isomorphism \cite{BC1}, \cite{BCH}. Here $ K_*(G \ltimes_\red A) $ is the $ K $-theory of the reduced crossed product of $ A $ by $ G $. 
The validity of this conjecture has applications in topology, geometry and representation theory. In particular, if $ G $ 
is discrete then the Baum-Connes conjecture with trivial coefficients $ \mathbb{C} $ 
implies the Novikov conjecture on higher signatures and the Kadison-Kaplansky idempotent 
conjecture. \\
Meyer and Nest have reformulated the Baum-Connes conjecture using the language
of triangulated categories and derived functors \cite{MNtriangulated}. 
In this approach the left hand side of the assembly map is identified with the 
localisation $ \LH F $ of the functor $ F(A) = K_*(G \ltimes_\red A) $ on the equivariant Kasparov category $ KK^G $. 
Among other things, this description allows to establish permanence properties of the conjecture in an efficient way. \\
In addition, the approach in \cite{MNtriangulated} is a natural starting point to study an analogue of the 
Baum-Connes conjecture for locally compact quantum groups. The usual definition of the left hand side of the conjecture 
is based on the universal space for proper actions, a concept which does not translate to the quantum setting 
in an obvious way. Following \cite{MNtriangulated}, one has to specify instead 
an appropriate subcategory of the equivariant Kasparov category corresponding to compactly induced actions 
in the group case. This approach has been implemented in \cite{MNcompact} where a strong form 
of the Baum-Connes conjecture for duals of compact groups is established. Duals of compact groups are, in a sense, the most basic
examples of discrete quantum groups. \\
In this paper we develop these ideas further and prove an analogue of the Baum-Connes conjecture for free orthogonal 
quantum groups. These quantum groups, introduced by Wang and van Daele \cite{Wang}, \cite{vDW}, can be considered as 
quantum analogs of orthogonal matrix Lie groups. If $ Q \in GL_n(\mathbb{C}) $ is a matrix satisfying $ Q \overline{Q} = \pm 1 $  
then the free orthogonal quantum group $ A_o(Q) $ is the universal $ C^* $-algebra generated by the 
entries of a unitary $ n \times n $-matrix $ u $ satisfying the relations $ u = Q \overline{u} Q^{-1} $. 
In this paper we will use the notation $ A_o(Q) = C^*_\max(\mathbb{F}O(Q)) $ in order to emphasize that we view 
this $ C^* $-algebra as the full group $ C^* $-algebra of a discrete quantum group. Accordingly, we will 
refer to $ \mathbb{F}O(Q) $ as the free orthogonal quantum group associated to $ Q $. 
In the case that $ Q = 1_n \in GL_n(\mathbb{C}) $ is the identity matrix we simply write $ \mathbb{F}O(n) $ instead 
of $ \mathbb{F}O(1_n) $. In fact, this special case illustrates the link to classical orthogonal groups since the algebra 
$ C(O(n)) $ of functions on $ O(n) $ is the abelianization of $ C^*_\max(\mathbb{F}O(n)) $. 
It is known \cite{Banicaunitary} that the quantum group $ \mathbb{F}O(Q) $ is not amenable if $ Q \in GL_n(\mathbb{C}) $ 
for $ n > 2 $. \\
The main result of this paper is that $ \mathbb{F}O(Q) $ has a $ \gamma $-element 
and that $ \gamma = 1 $ for $ Q \in GL_n(\mathbb{C}) $ and $ n > 2 $. The precise meaning of this statement, 
also referred to as the strong Baum-Connes conjecture, 
will be explained in section \ref{secbc} using the language of triangulated 
categories. However, we point out that triangulated categories are not needed to describe the applications 
that motivated our study. 
Firstly, it follows that free orthogonal quantum groups are $ K $-amenable. 
This answers in an affirmative way a question arising from the work of Vergnioux 
on quantum Cayley trees \cite{Vergniouxtrees}. 
Secondly, by studying the left hand side of the assembly map we obtain 
an explicit calculation of the $ K $-theory of $ \mathbb{F}O(Q) $. 
In the same way as in the case of free groups, the result of this calculation implies that 
the reduced group $ C^* $-algebra of $ \mathbb{F}O(n) $ does not contain nontrivial idempotents. 
This may be regarded as an analogue of the Kadison-Kaplansky conjecture. \\
Our results support the point of view that free quantum groups behave like free groups in many respects. 
By work of Vaes and Vergnioux \cite{VaesVergnioux} it is known, for instance, that the reduced $ C^* $-algebra 
$ C^*_\red(\mathbb{F}O(n)) $ of $ \mathbb{F}O(n) $ 
is exact and simple for $ n > 2 $. Moreover, the associated von Neumann algebra $ \mathcal{L}(\mathbb{F}O(n)) $ is a full and prime
factor. Let us note that, in contrast to the case of free groups, even the 
$ K $-theory of the maximal $ C^* $-algebras of orthogonal quantum groups seems difficult to compute directly. \\
As already mentioned above, our approach is motivated from the work of Meyer and Nest \cite{MNtriangulated}. 
In fact, the definition of the assembly map for torsion-free quantum groups proposed by Meyer in \cite{Meyerhomalg2} is 
the starting point of this paper. We proceed by observing that the strong Baum-Connes conjecture for 
torsion-free quantum groups is invariant under monoidal equivalence. 
The theory of monoidal equivalence for compact quantum groups was developed by Bichon, De Rijdt, and Vaes \cite{BdRV}. 
We use it to translate the Baum-Connes problem for free orthogonal quantum groups into a 
specific problem concerning $ SU_q(2) $. This step builds on the results in \cite{BdRV} and 
the foundational work of Banica \cite{Banicafo}.
The crucial part of our argument is a precise study of the equivariant $ KK $-theory of the standard Podle\'s sphere. 
By definition, the Podle\'s sphere $ SU_q(2)/T $ is the homogeneous space of $ SU_q(2) $ with respect to the classical 
maximal torus $ T \subset SU_q(2) $. Our constructions in connection with the Podle\'s sphere rely on the 
considerations in \cite{NVpoincare}. Finally, the $ K $-theory computation for $ \mathbb{F}O(Q) $ involves some
facts from homological algebra for triangulated categories worked out in \cite{Meyerhomalg2}. \\
Let us describe how the paper is organized. In section \ref{secspec} we discuss some preliminaries on 
compact quantum groups. In particular, we review the construction of 
spectral subspaces for actions of compact quantum groups on $ C^* $-algebras and Hilbert modules.
Section \ref{secsuq2} contains basic definitions related to $ SU_q(2) $ 
and the standard Podle\'s sphere $ SU_q(2)/T $. Moreover, it is shown that the dual of $ SU_q(2) $ can be viewed as 
a torsion-free discrete quantum group in the sense of \cite{Meyerhomalg2}. The most technical part of the paper is 
section \ref{secpodles} which contains our results on the Podle\'s sphere.  
In section \ref{secbc} we explain the formulation of the Baum-Connes conjecture 
for torsion-free quantum groups proposed in \cite{Meyerhomalg2}. Using the considerations from section 
\ref{secpodles} we prove that the dual of $ SU_q(2) $ satisfies the strong Baum-Connes conjecture in section \ref{secbcsuq2}. 
Section \ref{secfreeqg} contains the definition of free orthogonal quantum groups and a brief review of 
the theory of monoidal equivalence for compact quantum groups \cite{BdRV}. 
In section \ref{secmon} we show that monoidally equivalent compact quantum groups have equivalent 
equivariant $ KK $-categories. This implies that the strong Baum-Connes property is invariant under 
monoidal equivalence. Due to the work in \cite{BdRV} and our results in section \ref{secbcsuq2} it follows that 
free orthogonal quantum groups satisfy the strong Baum-Connes conjecture. Finally, in section \ref{secapp} 
we discuss applications and consequences. \\
Let us make some remarks on notation. We write $ \LH(\E) $ for the space of adjointable
operators on a Hilbert $ A $-module $ \E $. Moreover $ \KH(\E) $ denotes the
space of compact operators. The closed linear span of a subset $ X $ of a Banach space is denoted by $ [X] $. 
Depending on the context, the symbol
$ \otimes $ denotes either the tensor product of Hilbert spaces, the minimal tensor product of $ C^* $-algebras, or the tensor product
of von Neumann algebras. We write $ \odot $ for algebraic tensor products. For operators on multiple tensor products 
we use the leg numbering notation. \\
It is a pleasure to thank U. Krähmer, R. Meyer, R. Nest and N. Vander Vennet for helpful discussions on the subject of this paper. 

\section{Compact quantum groups and spectral decomposition} \label{secspec}

Concerning the general theory of quantum groups, we assume that the reader is familiar with the definitions and constructions 
that are reviewed in the first section of \cite{NVpoincare}. For more information we refer to the literature 
\cite{BSUM}, \cite{Kustermansuniversal}, \cite{KVLCQG}, \cite{Vaesimprimitivity}, \cite{Woronowiczleshouches}. Unless explicitly 
stated otherwise, our notation and conventions will follow \cite{NVpoincare} throughout the paper. \\
The purpose of this section is to explain some specific preliminaries on compact quantum groups. In 
particular, we discuss the construction of spectral subspaces for actions of compact quantum groups on $ C^* $-algebras and 
Hilbert modules. \\
Let $ G $ be a compact quantum group. Since $ G $ is compact the corresponding reduced $ C^* $-algebra of functions $ C^\red(G) $ 
is unital. A (unitary) representation $ \pi $ of $ G $ on a Hilbert space $ \H_\pi $ is an invertible 
(unitary) element $ u^\pi \in \LH(C^\red(G) \otimes \H_\pi) $ satisfying the relation
$$
(\Delta \otimes \id)(u^\pi) = u^\pi_{13} u^\pi_{23}.
$$
That is, a unitary representation of $ G $ is the same thing as a unitary corepresentation of the 
Hopf-$ C^* $-algebra $ C^\red(G) $. \\
Let $ \pi, \eta $ be representations of $ G $ on the Hilbert spaces
$ \H_\pi, \H_\eta $, given by the invertible elements $ u^\pi \in \LH(C^\red(G) \otimes \H_\pi) $ 
and $ u^\eta \in \LH(C^\red(G) \otimes \H_\eta) $, respectively. An operator $ T $ in $ \LH(\H_\pi, \H_\eta) $ is called an intertwiner 
between $ \pi $ and $ \eta $ if $ (\id \otimes T)u^\pi = u^\eta(\id \otimes T) $. We will denote the space of intertwiners 
between $ \H_\pi $ and $ \H_\eta $ by $ \Mor(\H_\pi, \H_\eta) $. The representations $ \pi $ and $ \eta $ are 
equivalent iff $ \Mor(\H_\pi, \H_\eta) $ contains an invertible operator. Every unitary representation of a compact quantum group 
decomposes into a direct sum of irreducibles, and all irreducible representations are finite dimensional. 
Every finite dimensional representation is equivalent to a 
unitary representation, and according to Schur's lemma a representation $ \pi $ 
is irreducible iff $ \dim(\Mor(\H_\pi, \H_\pi)) = 1 $. 
By slight abuse of notation, we will sometimes write $ \hat{G} $ for the set of 
isomorphism classes of irreducible unitary representations of $ G $. The trivial representation of $ G $ on the one-dimensional 
Hilbert space is denoted by $ \epsilon $. \\
The tensor product of the representations $ \pi $ and $ \eta $ is the representation 
$ \pi \otimes \eta $ on $ \H_\pi \otimes \H_\eta $ given 
by $ u^{\pi \otimes \eta} = u^\eta_{13} u^\pi_{12} \in \LH(C^\red(G) \otimes \H_\pi \otimes \H_\eta) $. 
The class of all finite dimensional representations of $ G $ together with the intertwining operators as morphisms and 
the direct sum and tensor product operations yields a $ C^* $-tensor category $ \Rep(G) $. 
This category is called the representation category of $ G $. By construction, it comes 
equipped with a canonical $ C^* $-tensor functor to the category of finite dimensional Hilbert spaces. \\
Let $ \pi $ be a finite dimensional representation given by $ u^\pi \in \LH(C^\red(G) \otimes \H_\pi) $, and let 
$ \dim(\pi) = \dim(\H_\pi) = n $ be the dimension of the underlying Hilbert space. 
If $ e^\pi_1, \dots, e^\pi_n $ is an orthonormal basis for $ \H_\pi $ we obtain associated matrix elements 
$ u^\pi_{ij} \in C^\red(G) $ given by 
$$
u^\pi_{ij} = \bra e^\pi_i, u^\pi(e^\pi_j) \ket.
$$ 
The corepresentation identity for $ u^\pi $ corresponds to 
$$ 
\Delta(u^\pi_{ij}) = \sum_{k = 1}^n u^\pi_{ik} \otimes u^\pi_{kj}
$$ 
for $ 1 \leq i, j \leq n $. Conversely, a (unitary) invertible matrix 
$ u = (u_{ij}) \in M_n(C^\red(G)) = \LH(C^\red(G) \otimes \mathbb{C}^n) $ 
satisfying these relations yields a (unitary) representation of $ G $. \\
The linear span of the matrix elements of $ \pi \in \hat{G} $ forms a finite dimensional coalgebra 
$ \mathbb{C}[G]_\pi \subset C^\red(G) $. Moreover 
$$ 
\mathbb{C}[G] = \bigoplus_{\pi \in \hat{G}} \mathbb{C}[G]_\pi
$$
is a dense Hopf $ * $-algebra $ \mathbb{C}[G] \subset C^\red(G) $ by the Peter-Weyl theorem. 
In subsequent sections we will use the fact that similar spectral decompositions exist for 
arbitrary $ G $-$ C^* $-algebras and $ G $-Hilbert modules. In order to discuss this we review some further definitions 
and results. \\
Let $ \pi $ be an irreducible unitary representation of $ G $, and let $ u^\pi_{ij} $ be the matrix elements in some fixed basis. 
The contragredient representation $ \pi^c $ is given by the matrix $ (u^{\pi^c})_{ij} = S(u^\pi_{ji}) $ where $ S $ is the antipode 
of $ \mathbb{C}[G] $. In general $ \pi^c $ is not unitary, but as any finite dimensional representation of $ G $ it is 
unitarizable. The representations $ \pi $ and $ \pi^{cc} $ are equivalent, and there exists a unique positive invertible 
intertwiner $ F_\pi \in \Mor(\H_\pi, \H_{\pi^{cc}}) $ satisfying $ \tr(F_\pi) = \tr(F_\pi^{-1}) $. The trace of $ F_\pi $ is 
called the quantum dimension of $ \pi $ and denoted by $ \dim_q(\pi) $. \\
With this notation, the Schur orthogonality relations are
$$
\phi(u^\pi_{ij} (u^\eta_{kl})^*) = \delta_{\pi\eta} \delta_{ik} \, \frac{1}{\dim_q(\pi)}\, (F_\pi)_{lj}
$$
where $ \pi, \eta \in \hat{G} $ and $ \phi: C^\red(G) \rightarrow \mathbb{C} $ is the Haar state of $ G $. 
In the sequel we shall fix bases such that $ F_\pi $ is a diagonal operator for all $ \pi \in \hat{G} $. \\
Let $ \pi $ be a unitary representation of $ G $ with matrix elements $ u^\pi_{ij} $. The element 
$$ 
\chi_\pi = \sum_{j = 1}^{\dim(\pi)} u^\pi_{jj}
$$
in $ C^\red(G) $ is called the character of $ \pi $. The subring in $ C^\red(G) $ generated by the 
characters of unitary representations is isomorphic to the (opposite of the) representation ring of $ G $. \\
Let us now fix our terminology concerning coactions. By an algebraic coaction of $ \mathbb{C}[G] $ on a vector space $ M $ we mean 
an injective linear map $ \gamma: M \rightarrow \mathbb{C}[G] \odot M $ such that $ (\id \odot \gamma) \gamma = (\Delta \odot \id) \gamma $. 
For an algebraic coaction one always has $ (\epsilon \odot \id) \gamma = \id $ and 
$ (\mathbb{C}[G] \odot 1)\gamma(M) = \mathbb{C}[G] \odot M $. 
Hence a vector space together with an algebraic coaction 
of $ \mathbb{C}[G] $ is the same thing as a (left) $ \mathbb{C}[G] $-comodule. 
Accordingly one defines right coactions on vector spaces. \\ 
By an algebraic coaction of $ \mathbb{C}[G] $ on a $ * $-algebra $ \A $ we shall mean an injective $ * $-homomorphism 
$ \alpha: \A \rightarrow \mathbb{C}[G] \odot \A $ such that $ (\id \odot \alpha) \alpha = (\Delta \odot \id)\alpha $. 
Accordingly one defines right coactions on $ * $-algebras. A $ * $-algebra equipped with an algebraic 
coaction of $ \mathbb{C}[G] $ will also be called a $ G $-algebra. \\
Let $ G $ be a compact quantum group and let $ A $ be a $ G $-$ C^* $-algebra 
with coaction $ \alpha: A \rightarrow M(C^\red(G) \otimes A) $. 
Since $ G $ is compact, such a coaction is an injective $ * $-homomorphism $ \alpha: A \rightarrow C^\red(G) \otimes A $ 
satisfying the coassociativity identity $ (\Delta \otimes \id) \alpha = (\id \otimes \alpha) \alpha $ 
and the density condition $ [(C^\red(G) \otimes 1)\alpha(A)] = C^\red(G) \otimes A $. For $ \pi \in \hat{G} $ we let 
$$
A_\pi = \{a \in A| \alpha(a) \in \mathbb{C}[G]_\pi \odot A \}
$$
be the $ \pi $-spectral subspace of $ A $. The subspace $ A_\pi $ is a closed in $ A $, and there is a 
projection operator $ p_\pi: A \rightarrow A_\pi $ defined by 
$$
p_\pi(a) = (\theta_\pi \otimes \id) \alpha(a)
$$
where 
$$ 
\theta_\pi(x) = \dim_q(\pi) \sum_{j = 1}^{\dim(\pi)} (F_\pi)^{-1}_{jj} \phi(x (u_{jj}^\pi)^*). 
$$
The spectral subalgebra $ \S(A) \subset A $ is the $ * $-subalgebra defined by 
$$
\S(A) = \S_G(A) = \bigoplus_{\pi \in \hat{G}} A_\pi, 
$$
and we note that $ \S(A) $ is a $ G $-algebra in a canonical way. 
From the Schur orthogonality relations and $ [(C^\red(G) \otimes 1)\alpha(A)] = C^\red(G) \otimes A $ 
it is easy to check that $ \S(A) $ is dense in $ A $, compare \cite{Podlessym}. \\
In a similar way one defines the spectral decomposition of $ G $-Hilbert modules. 
Let $ \E_A $ be a $ G $-Hilbert $ A $-module over the $ G $-$ C^* $-algebra $ A $. 
Since $ G $ is compact, the corresponding coaction $ \gamma: \E \rightarrow M(C^\red(G) \otimes \E) $ 
is an injective linear map $ \E \rightarrow C^\red(G) \otimes \E $ satisfying the coaction identity 
$ (\Delta \otimes \id) \gamma = (\id \otimes \gamma)\gamma $ and the density 
condition $ [(C^\red(G) \otimes 1)\gamma(\E)] = C^\red(G) \otimes \E $. For $ \pi \in \hat{G} $ we let 
$$
\E_\pi = \{\xi \in \E| \gamma(\xi) \in \mathbb{C}[G]_\pi \odot \E \} 
$$
be the corresponding spectral subspace. As in the algebra case, the spectral subspace $ \E_\pi $ is closed in 
$ \E $, and there is a projection map $ p_\pi: \E \rightarrow \E_\pi $. \\
By definition, the spectral submodule of $ \E $ is the dense subspace
$$
\S(\E) = \bigoplus_{\pi \in \hat{G}} \E_\pi
$$
of $ \E $. The spectral submodule $ \S(\E) $ is in fact a right $ \S(A) $-module, and the scalar product of $ \E $ restricts 
to an $ \S(A) $-valued inner product on $ \S(\E) $. In this way $ \S(\E) $ becomes a pre-Hilbert $ \S(A) $-module. \\
For a $ G $-algebra $ \A $ the spectral subspace $ \A_\pi $ for $ \pi \in \hat{G} $ is defined in the same way. 
This yields a corresponding spectral decomposition of $ \A $, the difference to the $ C^* $-setting is that we always 
have $ \S(\A) = \A $ in this case. The same remark applies to coactions of $ \mathbb{C}[G] $ on arbitrary vector spaces. 
If $ H $ is another compact quantum group and $ M $ is a $ \mathbb{C}[G] $-$ \mathbb{C}[H] $-bicomodule, we will also write 
$ _\pi M $ for the $ \pi $-spectral subspace corresponding to the left coaction. \\
Finally, we recall the definition of cotensor products. Let $ M $ be a right $ \mathbb{C}[G] $-comodule 
with coaction $ \rho: M \rightarrow M \odot \mathbb{C}[G] $ and $ N $ be a left $ \mathbb{C}[G] $-comodule 
with coaction $ \lambda: N \rightarrow \mathbb{C}[G] \odot N $. 
The cotensor product of $ M $ and $ N $ is the equalizer 
\begin{equation*}
    \xymatrix{
      {M \Box_{\mathbb{C}[G]} N} \ar@{->}[r] & {M \odot N \;} \ar@<1ex>@{->}[r] &
          {\;M \odot \mathbb{C}[G] \odot N} 
            \ar@<1ex>@{<-}[l]
               } 
\end{equation*}
of the maps $ \id \odot \lambda $ and $ \rho \odot \id $. 

\section{The quantum group $ SU_q(2) $} \label{secsuq2}

In this section we review some definitions and constructions related to $ SU_q(2) $ \cite{Woronowiczsuq2}. Moreover we 
show that the dual discrete quantum group of $ SU_q(2) $ is torsionfree in a suitable sense. 
Throughout we consider $ q \in [-1,1] \setminus \{0\} $, at some points we will exclude the cases $ q = \pm 1 $. \\
By definition, $ C(SU_q(2)) $ is the universal $ C^* $-algebra generated by elements $ \alpha $ and $ \gamma $ satisfying the relations
$$
\alpha \gamma = q \gamma \alpha, \quad \alpha \gamma^* = q \gamma^* \alpha, \quad \gamma \gamma^* = \gamma^* \gamma, \quad 
\alpha^* \alpha + \gamma^* \gamma = 1, \quad \alpha \alpha^*  + q^2 \gamma \gamma^* = 1. 
$$
These relations are equivalent to saying that the fundamental matrix 
$$
u = 
\begin{pmatrix}
\alpha & -q \gamma^* \\
\gamma & \alpha^*
\end{pmatrix}
$$
is unitary. \\
The comultiplication $ \Delta: C(SU_q(2)) \rightarrow C(SU_q(2)) \otimes C(SU_q(2)) $ is given on 
the generators by 
$$
\Delta(\alpha) = \alpha \otimes \alpha - q \gamma^* \otimes \gamma, \qquad \Delta(\gamma) = \gamma \otimes \alpha + \alpha^* \otimes \gamma, 
$$
and in this way the compact quantum group $ SU_q(2) $ is defined. We remark that 
there is no need to distinguish between full and reduced $ C^* $-algebras here since $ SU_q(2) $ is coamenable, see \cite{Banicafusion}. \\
The Hopf $ * $-algebra $ \mathbb{C}[SU_q(2)] $ is the dense $ * $-subalgebra of $ C(SU_q(2)) $ generated by $ \alpha $ and $ \gamma $, 
with counit $ \epsilon: \mathbb{C}[SU_q(2)] \rightarrow \mathbb{C} $ and antipode $ S: \mathbb{C}[SU_q(2)] \rightarrow \mathbb{C}[SU_q(2)] $ 
determined by
$$ 
\epsilon(\alpha) = 1, \qquad \epsilon(\gamma) = 0 
$$ 
and 
\begin{align*}
S(\alpha) = \alpha^*, \qquad S(\alpha^*) = \alpha, \qquad S(\gamma) = -q \gamma, \qquad S(\gamma^*) = -q^{-1} \gamma^*, 
\end{align*}
respectively.  
We use the Sweedler notation $ \Delta(x) = x_{(1)} \otimes x_{(2)} $ for the 
comultiplication and write 
$$ 
f \hit x = x_{(1)} f(x_{(2)}), \qquad x \hitby f = f(x_{(1)}) x_{(2)} 
$$ 
for elements $ x \in \mathbb{C}[SU_q(2)] $ and linear functionals $ f: \mathbb{C}[SU_q(2)] \rightarrow \mathbb{C} $. \\
The antipode is an algebra antihomomorphism satisfying 
$ S(S(x^*)^*) = x $ for all $ x \in \mathbb{C}[SU_q(2)] $. In particular the map $ S $ is invertible, and 
the inverse of $ S $ can be written as 
$$
S^{-1}(x) = f_1 \hit S(x) \hitby f_{-1}
$$
where $ f_1: \mathbb{C}[SU_q(2)] \rightarrow \mathbb{C} $ is the modular character given by 
$$
f_1(\alpha) = |q|^{-1}, \qquad f_1(\alpha^*) = |q|, \qquad f_1(\gamma) = 0, \qquad f_1(\gamma^*) = 0
$$
and $ f_{-1} $ is defined by $ f_{-1}(x) = f(S(x)) $. These maps are actually members of a canonical family 
$ (f_z)_{z \in \mathbb{C}} $ of characters. 
The character $ f_1 $ describes the modular properties of the Haar state $ \phi $ of $ C(SU_q(2)) $ in the sense that 
$$
\phi(xy) = \phi(y (f_1 \hit x \hitby f_1)) 
$$
for all $ x,y \in \mathbb{C}[SU_q(2)] $. \\
For $ q \in (-1,1) \setminus \{0\} $ we denote by $ U_q(\mathfrak{sl}(2))$ the quantized universal enveloping algebra of 
$ \mathfrak{sl}(2) $. This is the algebra generated by the elements $ E, F, K $ such that $ K $ is invertible and the relations 
\begin{align*}
K E K^{-1} = q^2 E, \qquad K F K^{-1} = q^{-2} F, \qquad [E,F] &= \frac{K - K^{-1}}{q - q^{-1}}
\end{align*}
are satisfied. We consider $ U_q(\mathfrak{sl}(2))$ with its Hopf algebra structure determined by
\begin{align*}
\Delta(K) = K \otimes K, \qquad \Delta (E) = E \otimes K &+ 1 \otimes E, \qquad \Delta(F) = F \otimes 1 + K^{-1} \otimes F \\ 
\epsilon(K) = 1, \qquad \epsilon(E) &= 0, \qquad \epsilon(F) = 0 \\ 
S(K) = K^{-1}, \qquad S(E) = & -E K^{-1}, \qquad S(F) = -K F
\end{align*}
and the $ * $-structure defining the compact real form, explicitly
$$
K^* = K, \qquad E^* = \sgn(q) FK, \qquad F^* = \sgn(q) K^{-1}E
$$
where $ \sgn(q) $ denotes the sign of $ q $. Let us remark that in the literature sometimes a wrong $ * $-structure 
is used in the case $ q < 0 $. \\
There is a nondegenerate skew-pairing between the Hopf-$ * $-algebras $ U_q(\mathfrak{sl}(2)) $ and 
$ \mathbb{C}[SU_q(2)] $, compare \cite{KS}. In particular, every finite dimensional unitary corepresentation of 
$ C(SU_q(2)) $ corresponds to a finite dimensional unital $ * $-representation of $ U_q(\mathfrak{sl}(2)) $. \\
Recall that a discrete group is called torsion-free if it does not contain nontrivial elements of finite order. 
For discrete quantum groups the following definition was proposed by Meyer \cite{Meyerhomalg2}. 
\begin{definition} \label{deftorsion}
A discrete quantum group $ G $ is called torsion-free iff every finite dimensional $ \hat{G} $-$ C^* $-algebra 
for the dual compact quantum group $ \hat{G} $ is isomorphic to a direct sum of $ \hat{G} $-$ C^* $-algebras that are 
equivariantly Morita equivalent to $ \mathbb{C} $. 
\end{definition}
In other words, according to definition \ref{deftorsion}, a discrete quantum group $ G $ is torsion-free if for every finite dimensional 
$ \hat{G} $-$ C^*$-algebra $ A $ there are finite dimensional Hilbert spaces $ \H_1, \dots, \H_l $ and 
unitary corepresentations $ u_j \in \LH(C^\red(\hat{G}) \otimes \H_j) $ such that $ A $  
is isomorphic to $ \KH(\H_1) \oplus \cdots \oplus \KH(\H_l) $ as a $ \hat{G} $-$ C^* $-algebra. Here 
each matrix block $ \KH(\H_j) $ is equipped with the adjoint action corresponding to 
$ u_j $. If $ G $ is a discrete group this is equivalent to the usual notion of torsion-freeness. \\ 
Definition \ref{deftorsion} is motivated from the study of torsion phenomena that occur for coactions of 
compact groups \cite{MNcompact}. Hence it is not surprising that it also provides the correct picture for duals of 
$ q $-deformations. We shall discuss explicitly the case of $ SU_q(2) $. 
\begin{prop} \label{suq2torsion}
Let $ q \in (-1,1) \setminus \{0\} $. Then the discrete quantum group dual to $ SU_q(2) $ is torsion-free. 
\end{prop}
\proof The following argument was suggested by U. Kr\"ahmer. For simplicity we restrict ourselves to the case $ q > 0 $, 
the case of negative $ q $ is treated in a similar way. Let us assume that $ A $ is a finite dimensional $ SU_q(2) $-$ C^* $-algebra 
with coaction $ \alpha: A \rightarrow C(SU_q(2)) \otimes A $. Since $ A $ is finite dimensional we may 
write $ A = M_{n_1}(\mathbb{C}) \oplus \cdots \oplus M_{n_l}(\mathbb{C}) $. The task is 
to describe the coaction in terms of this decomposition. \\
First consider the restriction of $ \alpha $ to $ C(T) $. Since the torus
$ T $ is a connected group, the corresponding action of $ T $ preserves the decomposition of $ A $ into matrix blocks. Moreover, 
on each matrix block the action arises from a representation of $ T $. 
Hence, if we view $ A $ as a $ U_q(\mathfrak{sl}(2)) $-module algebra, the action of $ K $ is implemented by conjugating with an 
invertible self-adjoint matrix $ k = k_1 \oplus \cdots \oplus k_l $. Moreover we may suppose that $ k $ has only positive eigenvalues. With 
these requirements each of the matrices $ k_j $ is uniquely determined up to a positive scalar factor. \\
Next consider the skew-primitive elements $ E $ and $ F $. From the definition of the comultiplication in $ U_q(\mathfrak{sl}(2)) $ we obtain
$$
E \cdot (ab) = (E \cdot a)(K \cdot b) + a(E \cdot b), \qquad F \cdot(ab) = F(a) b + (K^{-1} \cdot a) (F \cdot b)
$$
for all $ a, b \in A $. Hence $ E $ and $ F $ can be viewed as Hochschild-$ 1 $-cocycles on $ A $ 
with values in appropriate $ A $-bimodules. Since $ A $ is a semisimple algebra the corresponding Hochschild cohomology groups vanish. 
Consequently there are $ e, f \in A $ such that 
\begin{align*} 
E \cdot a &= e k^{-1} (K \cdot a) - a e k^{-1}\\
F \cdot a &= f a - (K^{-1} \cdot a) f
\end{align*}
for all $ a \in A $. 
It follows that $ E $ and $ F $ preserve the decomposition of $ A $ into matrix blocks. In particular, we may restrict attention 
to the case that $ A $ is a simple matrix algebra. \\
Let us assume $ A = M_n(\mathbb{C}) $ in the sequel. Then $ e $ and $ f $ are uniquely determined up to addition of a scalar multiple of $ 1 $ and $ k^{-1} $, respectively. 
The relation $ K E K^{-1} = q^2 E $ implies
\begin{align*}
k(e k^{-1} a - k^{-1} a k ek^{-1}) k^{-1} &= k e k^{-1} a k^{-1} - a k e k^{-2} = q^2 (e a k^{-1} - a e k^{-1})  
\end{align*}
for all $ a \in A $ and yields $ k e k^{-1} - q^2 e = \lambda $ 
for some $ \lambda \in \mathbb{C} $. If we replace $ e $ by $ e - \lambda (1- q^2)^{-1} $ we obtain 
$ k e k^{-1} = q^2 e $. Similarly we may achieve $ k f k^{-1} = q^{-2} f $, and we fix 
$ e $ and $ f $ such that these identities hold. The commutation relation for $ E $ and $ F $ 
implies
\begin{align*}
\biggl(e&(fa - k^{-1} a  kf)k^{-1} - (fa - k^{-1} a  kf)e k^{-1}\biggr) \\
&\qquad - \biggl(f(e a k^{-1} - a e k^{-1}) - k^{-1}(e a k^{-1} - a e k^{-1}) kf\biggr) \\
&= (ef - fe) a k^{-1} - k^{-1} a (ef - fe) \\
&= \frac{1}{q - q^{-1}}\biggl(k a k^{-1} - k^{-1} a k\biggr).
\end{align*}
As a consequence we obtain
$$
e f - f e - \frac{k}{q - q^{-1}}= - \frac{\mu}{q - q^{-1}} k^{-1}
$$
for some constant $ \mu $. In fact, since $ k $ has only positive eigenvalues the scalar $ \mu $ is strictly positive. 
Replacing $ k $ by $ \lambda k $ and $ e $ by $ \lambda e $ with $ \lambda = \mu^{-1/2} $ 
yields
$$
[e,f] = \frac{k - k^{-1}}{q - q^{-1}}.
$$
It follows that there is a representation of $ U_q(\mathfrak{sl}(2)) $ on $ \mathbb{C}^n $ 
which induces the given action on $ A $ by conjugation. \\
We have $ (E \cdot a)^* = - F \cdot a^* $ for all $ a \in A $ and hence
$$
a^* (e^* - kf)  = (e^* - kf) a^*
$$
which implies $ e^* - kf = \nu $ for some $ \nu \in \mathbb{C} $. 
Conjugating with $ k $ yields
$$
\nu = k^{-1}(e^* - kf)k = (k e k^{-1})^* - q^2 kf = q^2 (e^* - kf) = q^2 \nu
$$
and thus $ \nu = 0 $. It follows that the representation of $ U_q(\mathfrak{sl}(2)) $ 
given by $ e, f $ and $ k $ is a $ * $-representation. We conclude that there exists a unitary corepresentation of 
$ C(SU_q(2)) $ on $ \mathbb{C}^n $ which implements the coaction on $ A $ as desired. \qed \\
Let us next discuss the regular representation of $ SU_q(2) $ for $ q \in (-1,1) \setminus \{0\} $. 
We write $ L^2(SU_q(2)) $ for the Hilbert space obtained using the inner product 
$$ 
\bra x, y \ket = \phi(x^* y) 
$$ 
on $ C(SU_q(2)) $. By definition, the regular representation on $ L^2(SU_q(2)) $ is given by the multiplicative unitary 
$ W \in \LH(C(SU_q(2)) \otimes L^2(SU_q(2)) $. \\
The Peter-Weyl theory describes the decomposition of this representation into irreducibles. 
As in the classical case, the irreducible representations of $ SU_q(2) $ are labelled by half-integers $ l $, and 
the corresponding Hilbert spaces have dimension $ 2l + 1 $. 
The matrix elements $ u^{(l)}_{ij} $ with respect to weight bases determine an 
orthogonal set in $ L^2(SU_q(2)) $. 
Moreover, if we write 
$$
[a] = \frac{q^a - q^{-a}}{q - q^{-1}}
$$ 
for the $ q $-number associated to $ a \in \mathbb{C} $, then the vectors
$$
e^{(l)}_{i,j} = q^i [2l + 1]^{\frac{1}{2}} u^{(l)}_{i,j} 
$$
form an orthonormal basis of $ L^2(SU_q(2)) $, compare \cite{KS}. \\
The regular representation of $ C(SU_q(2)) $ on $ L^2(SU_q(2)) $ is given by 
\begin{align*}
\alpha \, e_{i,j}^{(l)} &= a_+(l,i,j) \, e_{i - \frac{1}{2}, j - \frac{1}{2}}^{\left(l + \frac{1}{2} \right)} 
+ a_- (l,i,j) \, e_{i - \frac{1}{2}, j - \frac{1}{2}}^{\left(l - \frac{1}{2} \right)} \\
\gamma \, e_{i,j}^{(l)} &= c_+ (l,i,j) \, e_{i + \frac{1}{2}, j - \frac{1}{2}}^{\left(l + \frac{1}{2} \right)} + 
c_- (l,i,j) \, e_{i + \frac{1}{2}, j - \frac{1}{2}}^{\left(l - \frac{1}{2} \right)}
\end{align*}
where the explicit form of $ a_\pm $ and $ c_\pm $ is
\begin{align*}
a_+(l,i,j) &= q^{2l + i + j + 1} \, 
\frac{(1 - q^{2l - 2j + 2})^{\tfrac{1}{2}} (1 - q^{2l - 2i + 2})^{\tfrac{1}{2}}}{ (1 - q^{4l + 2})^{\tfrac{1}{2}} (1 - q^{4l + 4})^{\tfrac{1}{2}}} \\
a_-(l,i,j) &= \frac{(1 - q^{2l + 2j})^{\tfrac{1}{2}} (1 - q^{2l + 2i})^{\tfrac{1}{2}}}{(1 - q^{4l})^{\tfrac{1}{2}} (1 - q^{4l+2})^{\tfrac{1}{2}}}
\end{align*}
and
\begin{align*}
c_+(l,i,j) &= - q^{l + j} \, \frac{(1 - q^{2l - 2j + 2})^{\tfrac{1}{2}} 
(1 - q^{2l + 2i + 2})^{\tfrac{1}{2}}}{(1 - q^{4l + 2})^{\tfrac{1}{2}} (1 - q^{4l + 4})^{\tfrac{1}{2}}} \\
c_-(l,i,j) &= q^{l + i} \, \frac{(1 - q^{2l + 2j})^{\tfrac{1}{2}} 
(1 - q^{2l - 2i})^{\tfrac{1}{2}}}{(1 - q^{4l})^{\tfrac{1}{2}} (1 - q^{4l + 2})^{\tfrac{1}{2}}}.
\end{align*}
Note here that $ a_-(l,i,j) $ vanishes if $ i = -l $ or $ j = -l $, and similarly, $ c_-(l,i,j) $ vanishes 
for $ i = l $ or $ j = -l $. We obtain
\begin{align*}
\alpha^* \, e_{i,j}^{(l)} = a_+(l - \tfrac{1}{2}, & i + \tfrac{1}{2}, j + \tfrac{1}{2})\; 
e^{\left(l - \frac{1}{2}\right)}_{i + \frac{1}{2}, j + \frac{1}{2}} 
+ a_-(l + \tfrac{1}{2},i + \tfrac{1}{2}, j + \tfrac{1}{2})\; e^{\left(l + \frac{1}{2}\right)}_{i + \frac{1}{2}, j + \frac{1}{2}}
\end{align*}
and 
\begin{align*}
\gamma^* \, e_{i,j}^{(l)} = c_+(l - \tfrac{1}{2}, & i - \tfrac{1}{2}, j + \tfrac{1}{2})\; 
e^{\left(l - \frac{1}{2}\right)}_{i - \frac{1}{2}, j + \frac{1}{2}} 
+ c_-(l + \tfrac{1}{2},i - \tfrac{1}{2}, j + \tfrac{1}{2})\; e^{\left(l + \frac{1}{2}\right)}_{i - \frac{1}{2}, j + \frac{1}{2}}
\end{align*}
for the action of $ \alpha^* $ and $ \gamma^* $, respectively. Let us also record the formulas 
$$
u^{(l)*}_{i,j} = (-1)^{2l + i + j} q^{j - i} u^{(l)}_{-i, - j}
$$
and
$$
e^{(l)*}_{i,j} = (-1)^{2l + i + j} q^{i + j} e^{(l)}_{-i, - j}
$$
for the adjoint. \\
The classical torus $ T = S^1 $ is a closed quantum subgroup of $ SU_q(2) $.
Explicitly, the inclusion $ T \subset SU_q(2) $ is determined by the $ * $-homomorphism
$ \pi: \mathbb{C}[SU_q(2)] \rightarrow \mathbb{C}[T] = \mathbb{C}[z,z^{-1}] $ given by
$$
\pi(\alpha) = z, \qquad \pi(\gamma) = 0. 
$$
By definition, the standard Podle\'s sphere $ SU_q(2)/T  $ is the corresponding homogeneous space \cite{Podlesspheres}.
In the algebraic setting, it is described by the dense
$ * $-subalgebra $ \mathbb{C}[SU_q(2)/T] \subset C(SU_q(2)/T) $ of coinvariants in $ \mathbb{C}[SU_q(2)] $
with respect to the right coaction $ (\id \otimes \pi)\Delta $ of $ \mathbb{C}[T] $. 
Equivalently, it is the unital $ * $-subalgebra of $ \mathbb{C}[SU_q(2)] $ generated by the elements 
$ A = \gamma^* \gamma $ and $ B = \alpha^* \gamma $. 
These elements satisfy the relations 
$$
A = A^*, \quad AB = q^2 BA, \quad BB^* = q^{-2}A(1 - A), \quad B^*B = A(1 - q^2 A), 
$$
and we record the following explicit formulas for the action of $ A $ and $ B $ on $ L^2(SU_q(2)) $ 
for $ q \in (-1,1) \setminus \{0\} $. Firstly, 
\begin{align*}
&\gamma^* \gamma\, e_{i,j}^{(l)} = \gamma^* \biggl(c_+(l,i,j)\; e^{\left(l + \frac{1}{2}\right)}_{i + \frac{1}{2}, j - \frac{1}{2}} + 
c_-(l,i,j) \; e^{\left(l - \frac{1}{2}\right)}_{i + \frac{1}{2}, j - \frac{1}{2}} \biggr) \\
&= c_+(l - 1, i, j) c_-(l,i,j) \; e^{(l - 1)}_{ij} + (c_+(l,i,j)^2 + c_-(l,i,j)^2)\; e^{(l)}_{ij} \\
&\qquad + c_-(l + 1, i, j) c_+(l,i,j) \; e^{(l + 1)}_{ij} \\
&= -q^{2l + i + j - 1} \, \frac{(1 - q^{2l - 2j})^{\tfrac{1}{2}} (1 - q^{2l + 2i})^{\tfrac{1}{2}}(1 - q^{2l + 2j})^{\tfrac{1}{2}} (1 - q^{2l - 2i})^{\tfrac{1}{2}}}
{(1 - q^{4l - 2})^{\tfrac{1}{2}} (1 - q^{4l})(1 - q^{4l + 2})^{\tfrac{1}{2}}}\, e^{(l - 1)}_{ij} \\
&\quad + \biggl(q^{2l + 2j} \, \frac{(1 - q^{2l - 2j + 2})(1 - q^{2l + 2i + 2})}{(1 - q^{4l + 2}) (1 - q^{4l + 4})} 
+ q^{2l + 2i} \, \frac{(1 - q^{2l + 2j})(1 - q^{2l - 2i})}{(1 - q^{4l}) (1 - q^{4l + 2})} \biggr) e^{(l)}_{ij} \\
&\;\; - q^{2l + i + j + 1} \, \frac{(1 - q^{2l + 2j + 2})^{\tfrac{1}{2}} 
(1 - q^{2l - 2i + 2})^{\tfrac{1}{2}}(1 - q^{2l - 2j + 2})^{\tfrac{1}{2}} 
(1 - q^{2l + 2i + 2})^{\tfrac{1}{2}}}{(1 - q^{4l + 4})(1 - q^{4l + 6})^{\tfrac{1}{2}}(1 - q^{4l + 2})^{\tfrac{1}{2}}}\, e^{(l + 1)}_{ij} 
\end{align*}
determines the action of $ A $. 
Similarly, we get
\begin{align*}
&\alpha^* \gamma\, e_{i,j}^{(l)} = \alpha^* \biggl(c_+(l,i,j) e^{\left(l + \frac{1}{2}\right)}_{i + \frac{1}{2}, j - \frac{1}{2}} 
+ c_-(l,i,j) e^{\left(l - \frac{1}{2}\right)}_{i + \frac{1}{2}, j - \frac{1}{2}} \biggr) \\
&= a_+(l - 1, i + 1, j) c_-(l,i,j) e^{(l - 1)}_{i + 1,j} \\
&\qquad + \biggl(a_+(l,i + 1,j) c_+(l,i,j) + a_-(l,i + 1,j) c_-(l,i,j) \biggr) e^{(l)}_{i + 1, j} \\
&\qquad + a_-(l + 1, i + 1, j) c_+(l,i,j) e^{(l + 1)}_{i + 1, j} \\
&= q^{3l + 2i + j} \, \frac{(1 - q^{2l - 2j})^{\tfrac{1}{2}} (1 - q^{2l - 2i - 2})^{\tfrac{1}{2}} (1 - q^{2l + 2j})^{\tfrac{1}{2}} 
(1 - q^{2l - 2i})^{\tfrac{1}{2}}}
{(1 - q^{4l - 2})^{\tfrac{1}{2}} (1 - q^{4l}) (1 - q^{4l + 2})^{\tfrac{1}{2}}} e^{(l - 1)}_{i + 1,j} \\
&\quad + 
\biggl(q^{l + i} \frac{(1 - q^{2l + 2j}) (1 - q^{2l + 2i + 2})^{\tfrac{1}{2}}(1 - q^{2l - 2i})^{\tfrac{1}{2}}}{(1 - q^{4l})(1 - q^{4l+2})} \\
&\qquad\quad -q^{3l + i + 2j + 2} \, \frac{(1 - q^{2l - 2j + 2}) (1 - q^{2l - 2i})^{\tfrac{1}{2}}(1 - q^{2l + 2i + 2})^{\tfrac{1}{2}}}
{(1 - q^{4l + 2})(1 - q^{4l + 4})} \biggr) e^{(l)}_{i + 1,j}  \\
&\quad - q^{l + j} \frac{(1 - q^{2l + 2j + 2})^{\tfrac{1}{2}} (1 - q^{2l + 2i + 4})^{\tfrac{1}{2}} (1 - q^{2l - 2j + 2})^{\tfrac{1}{2}} 
(1 - q^{2l + 2i + 2})^{\tfrac{1}{2}}}{(1 - q^{4l + 4})(1 - q^{4l + 6})^{\tfrac{1}{2}} (1 - q^{4l + 2})^{\tfrac{1}{2}}}  e^{(l + 1)}_{i + 1, j} 
\end{align*}
for the action of $ B $. \\
In the sequel we abbreviate $ SU_q(2) = G_q $. If $\mathbb{C}_k $ denotes the irreducible representation of $ T $ of weight $ k \in \mathbb{Z} $, 
then the cotensor product
\begin{equation*}
\Gamma(\E_k) = \Gamma(G_q \times_T \mathbb{C}_k) = \mathbb{C}[G_q] \Box_{\mathbb{C}[T]} \mathbb{C}_k \subset \mathbb{C}[G_q] 
\end{equation*}
is a noncommutative analogue of the space of sections of the homogeneous vector bundle $ G \times_T \mathbb{C}_k $ over $ G/T $.
The space $ \Gamma(\E_k) $ is a $ \mathbb{C}[G_q/T] $-bimodule in a natural way which is finitely 
generated and projective
both as a left and right $ \mathbb{C}[G_q/T]$-module. The latter follows from the fact that
$ \mathbb{C}[G_q/T] \subset \mathbb{C}[G_q] $ is a faithfully flat Hopf-Galois extension, compare \cite{Schauenburg}. \\
We denote by $ C(\E_k) $ the closure of $ \Gamma(\E_k) $ inside $ C(G_q) $. 
The space $ C(\E_k) $ is a $ G_q $-equivariant Hilbert $ C(G_q/T) $-module with coaction induced by comultiplication.
We write $ L^2(\E_k) $ for the $ G_q $-Hilbert space obtained by taking the closure of $ \Gamma(\E_k) $ 
inside $ L^2(G_q) $. \\
Let us recall the definition of the Drinfeld double $ \DD(G_q) $ of $ G_q $. It is the locally 
compact quantum group determined by $ C_0(\DD(G_q)) = C(G_q) \otimes C^*(G_q) $ with comultiplication
$$
\Delta_{\DD(G_q)} = (\id \otimes \sigma \otimes \id)(\id \otimes \ad(W) \otimes \id)(\Delta \otimes \hat{\Delta}). 
$$
Here $ \ad(W) $ denotes the adjoint action of the left regular multiplicative unitary $ W \in M(C(G_q) \otimes C^*(G_q)) $ 
and $ \sigma $ is the flip map. \\
Observe that both $ C(G_q) $ and $ C^*(G_q) $ are quotient Hopf-$ C^* $-algebras of $ C_0(\DD(G_q)) $. 
It is shown in \cite{NVpoincare} that a $ \DD(G_q) $-$ C^* $-algebra $ A $ is uniquely determined by coactions 
$ \alpha: A \rightarrow M(C(G_q) \otimes A) $ and $ \lambda: A \rightarrow M(C^*(G_q) \otimes A) $ 
satisfying the Yetter-Drinfeld compatibility condition
$$
(\sigma \otimes \id)(\id \otimes \alpha)\lambda = (\ad(W) \otimes \id) (\id \otimes \lambda)\alpha.
$$
In a similar fashion on can study $ \DD(G_q) $-equivariant Hilbert modules. 
For instance, the space $ C(\E_k) $ defined above carries a coaction $ \lambda: C(\E_k) \rightarrow M(C^*(G_q) \otimes C(\E_k)) $
given by $ \lambda(f) = \hat{W}^*(1 \otimes f) \hat{W} $ where $ \hat{W} = \Sigma W^* \Sigma $. 
Together with the canonical coaction of $ C(G_q) $ this turns $ C(\E_k) $ into a $ \DD(G_q) $-equivariant Hilbert $ C(G_q/T) $-module. 

\section{Equivariant $ KK $-theory for the Podle\'s sphere} \label{secpodles}

In this section we study the equivariant $ KK $-theory of the standard Podle\'s sphere $ SU_q(2)/T $. 
Background information on equivariant $ KK $-theory for quantum group actions can be found in 
\cite{BSKK}, \cite{NVpoincare}. Most of the ingredients in this section have already been introduced in 
\cite{NVpoincare} in the case $ q \in (0,1) $. However, for the purposes of this paper we have to allow for 
negative values of $ q $. In the sequel we consider $ q \in (-1,1) \setminus \{0\} $, and 
as in the previous section we abbreviate $ G_q = SU_q(2) $. \\
Let us first recall the definition of the Fredholm module corresponding to the Dirac operator on the standard 
Podle\'s sphere $ G_q/T $, compare \cite{DSPodles}, \cite{NVpoincare}. The underlying graded $ G_q $-Hilbert space 
is 
$$ 
\H = \H_1 \oplus \H_{-1} = L^2(\E_1) \oplus L^2(\E_{-1})
$$ 
with its natural coaction of $ C(G_q) $. The covariant representation $ \phi $ of 
$ C(G_q/T) $ is given by left multiplication. \\
We note that $ \H_1 $ and $ \H_{-1} $ are isomorphic representations of 
$ G_q $ due to Frobenius reciprocity. Hence we obtain a self-adjoint unitary operator $ F $ on $ \H $ by 
$$
F =
\begin{pmatrix}
0 & 1 \\
1 & 0
\end{pmatrix}
$$
if we identify the basis vectors $ e^{(l)}_{i,1/2} $ and $ e^{(l)}_{i,- 1/2} $ in $ \H_1 $ and $ \H_{-1} $, respectively. 
Using the explicit formulas for the generators of the Podle\'s sphere from section \ref{secsuq2} one checks that 
$ D = (\H, \phi, F) $ is a $ G_q $-equivariant Fredholm module defining an element 
in $ KK^{G_q}(C(G_q/T), \mathbb{C}) $. \\
Our first aim is to lift this construction to $ \DD(G_q) $-equivariant $ KK $-theory where $ \DD(G_q) $ denotes the Drinfeld double 
of $ G_q $ as in section \ref{secsuq2}. 
We recall that the $ C^* $-algebra  $ C(G_q/T) = \ind_T^{G_q}(\mathbb{C}) $ is naturally a $ \DD(G_q) $-$ C^* $-algebra \cite{NVpoincare}. 
The corresponding coaction $ C(G_q/T) \rightarrow M(C^*(G_q) \otimes C(G_q/T)) $ is determined by the adjoint action 
$$
h \cdot g = h_{(1)} g S(h_{(2)}) 
$$
of $ \mathbb{C}[G_q] $ on $ \mathbb{C}[G_q/T] $. \\
The following lemma shows that the Hilbert spaces $ L^2(\E_k) $ become $ \DD(G_q) $-Hilbert spaces in a similar way. 
\begin{lemma}\label{lemma0}
For every $ k \in \mathbb{Z} $ the formula 
$$
\omega(x)(\xi)= x_{(1)} \xi\, f_1 \hit S(x_{(2)})
$$
determines a $ * $-homomorphism $ \omega: C(G_q) \rightarrow \LH(L^2(\E_k)) $ which
turns $ L^2(\E_k) $ into a $ \DD(G_q) $-Hilbert space such that the representation of $ C(G_q/T) $ 
by left multiplication is covariant. 
\end{lemma}
\proof Let us write $ \H = L^2(\E_k) $ and define a map $ \omega: \mathbb{C}[G_q] \rightarrow \LH(\H) $ 
by the above formula, where we recall that $ f_1: \mathbb{C}[G_q] \rightarrow \mathbb{C} $ denotes the modular character. 
For $ \xi \in \Gamma(\E_k) $ one obtains
\begin{align*}
(\id \otimes \pi)\Delta(\omega(x)(\xi)) &= (\id \otimes \pi)\Delta(x_{(1)} \xi\, f_1 \hit S(x_{(2)})) \\
&= x_{(1)} \xi_{(1)} S(x_{(4)}) \otimes \pi(x_{(2)} \xi_{(2)} f_1 \hit S(x_{(3)})) \\
&= x_{(1)} \xi_{(1)} S(x_{(4)}) \otimes \pi(x_{(2)}) \pi(\xi_{(2)}) \pi(f_1 \hit S(x_{(3)})) \\
&= x_{(1)} \xi S(x_{(4)}) \otimes \pi(x_{(2)} f_1 \hit S(x_{(3)})) z^k \\
&= x_{(1)} \xi S(x_{(4)}) \otimes \pi(S^{-1}(x_{(3)}) \hitby f_1) \pi(x_{(2)}) z^k \\
&= x_{(1)} \xi S(x_{(3)}) \otimes \pi(f_{-1}(x_{(2)}) 1) z^k \\
&= x_{(1)} \xi f_1 \hit S(x_{(2)}) \otimes z^k \\
&= \omega(x)(\xi) \otimes z^k
\end{align*}
using that $ \mathbb{C}[T] $ is commutative 
and $ S^{-1}(x) = f_1 \hit S(x) \hitby f_{-1} $ for all $ x \in \mathbb{C}[G_q] $. 
It follows that the map $ \omega $ is well-defined. Using $ f_1 \hit x^* = (f_{-1} \hit x)^* $ and the 
modular properties of the Haar state $ \phi $ we obtain
\begin{align*}
\bra \omega(x^*)(\xi), \eta \ket &= \phi((\omega(x^*)(\xi))^* \eta) \\
&= \phi((x_{(1)}^* \xi f_1 \hit S(x_{(2)}^*))^* \eta) \\
&= \phi((x_{(1)}^* \xi f_1 \hit S^{-1}(x_{(2)})^*)^* \eta) \\
&= \phi(f_{-1} \hit S^{-1}(x_{(2)}) \xi^* x_{(1)} \eta) \\
&= \phi(S(x_{(2)}) \hitby f_{-1} \xi^* x_{(1)} \eta) \\
&= \phi(\xi^* x_{(1)} \eta\, f_1 \hit S(x_{(2)})) \\
&= \phi(\xi^* \omega(x)(\eta)) \\
&= \bra \xi, \omega(x)(\eta) \ket, 
\end{align*}
and this shows that $ \omega $ defines a $ * $-homomorphism from $ C(G_q) $ to $ \LH(\H) $. 
We deduce that $ \omega $ corresponds to a coaction 
$ \lambda: \H \rightarrow M(C^*(G_q) \otimes \H) $, and
$ \lambda $ combines with the standard coaction of $ C(G_q) $ such that $ \H $ becomes a $ \DD(G_q) $-Hilbert space. \\
Recall that the action of $ C(G_q/T) $ on $ \H $ by left multiplication yields a $ G_q $-equivariant 
$ * $-homomorphism $ \phi: C(G_q/T) \rightarrow \LH(\H) $. We have 
\begin{align*}
\omega(x)(\phi(g)(\xi)) &= x_{(1)} g\xi f_1 \hit S(x_{(2)}) \\
&= x_{(1)} g S(x_{(2)}) x_{(3)} \xi f_1 \hit S(x_{(4)}) \\
&= \phi(x_{(1)} \cdot g)(\omega(x_{(2)})(\xi))
\end{align*}
for all $ x \in \mathbb{C}[G_q], g \in \mathbb{C}[G_q/T] $ and $ \xi \in \H $, and this implies that $ \phi $ 
is covariant with respect to $ \lambda $. \qed \\
We shall now show that the Fredholm module constructed in the beginning of this section 
determines a $ \DD(G_q) $-equivariant $ KK $-element. 
\begin{prop} \label{DiracYD}
Let $ q \in (-1,1) \setminus \{0\} $. The Fredholm module $ D $ defined above induces an element $ [D] $ in 
$ KK^{\DD(G_q)}(C(G_q/T), \mathbb{C}) $ in a natural way. 
\end{prop}
\proof We have to verify that $ F $ commutes with the action of $ \DD(G_q) $
up to compact operators. Since $ F $ is $ G_q $-equivariant this amounts to 
showing 
$$
(C^*(G_q) \otimes 1)(1 \otimes F - \ad_\lambda(F)) \subset C^*(G_q) \otimes \KH(\H)
$$
where $ \lambda: \H \rightarrow M(C^*(G_q) \otimes \H) $ is the coaction on $ \H = \H_1 \oplus \H_{-1} $ 
defined in lemma \ref{lemma0}. It suffices to check that $ F $ commutes with the corresponding action 
$ \omega: C(G_q) \rightarrow \LH(\H) $ up to compact operators. 
This is an explicit calculation using the formulas for the regular representation from section \ref{secsuq2}. 
In fact, we will obtain the assertion as a special case of our computations below. It turns out that 
$ F $ is actually $ \DD(G_q) $-equivariant. \qed \\
In the sequel we need variants of the representations defined in lemma \ref{lemma0}. 
Let $ t \in [0,1] $ and consider the representation $ \pi_t $ of $ \mathbb{C}[G_q] $ 
on $ \Gamma(\E_k) $ given by 
$$
\pi_t(x)(\xi) = x_{(1)} \xi\, f_t \hit S(x_{(2)}) 
$$
where $ f_t: \mathbb{C}[G_q] \rightarrow \mathbb{C} $ is the modular character given by $ f_t(\alpha) = |q|^t \alpha $ and $ f_t(\gamma) = 0 $. 
This action has the correct algebraic properties to 
turn $ L^2(\E_k) $ into a $ \DD(G_q) $-Hilbert space except that it is not compatible with 
the $ * $-structures for $ t < 1 $. In order to proceed we need explicit formulas for the action of the generators. \\
More precisely, a straightforward computation based on the formulas for the regular representation in section \ref{secsuq2} yields 
\begin{align*}
&\pi_t(\alpha)(e^{(l)}_{i,j}) = |q|^t \alpha e_{i,j}^{(l)} \alpha^* + |q|^{-t} q^2 \gamma^* e_{i,j}^{(l)} \gamma \\
&= |q|^{t} q^{-1} a_+(l,-i,-j) \biggl(a_+(l + \tfrac{1}{2},i + \tfrac{1}{2},j + \tfrac{1}{2}) e_{i, j}^{\left(l + 1\right)} 
+ a_-(l + \tfrac{1}{2},i + \tfrac{1}{2},j + \tfrac{1}{2}) e_{i, j}^{\left(l\right)} \biggr) \\
& + |q|^{t} q^{-1} a_-(l,-i,-j) \biggl(a_+(l - \tfrac{1}{2},i + \tfrac{1}{2},j + \tfrac{1}{2}) e_{i, j}^{\left(l\right)} + 
a_-(l - \tfrac{1}{2},i + \tfrac{1}{2},j + \tfrac{1}{2}) e_{i, j}^{\left(l - 1 \right)} \biggr) \allowdisplaybreaks[2]\\
& - |q|^{-t} q^2 c_+(l - \tfrac{1}{2}, -i - \tfrac{1}{2}, -j + \tfrac{1}{2}) 
\biggl(c_+(l - 1, i, j) e^{\left(l - 1 \right)}_{i, j} + c_-(l,i, j) e^{\left(l\right)}_{i, j} \biggr)\\
& - |q|^{-t} q^2 c_-(l + \tfrac{1}{2}, -i - \tfrac{1}{2}, -j + \tfrac{1}{2}) 
\biggl(c_+(l, i, j)  e^{\left(l\right)}_{i, j} + c_-(l + 1,i , j)\; e^{\left(l + 1\right)}_{i, j} \biggr)
\end{align*} 
and
\begin{align*}
&\pi_t(\alpha^*)(e^{(l)}_{i,j}) = |q|^{-t} \alpha^* e_{ij}^{(l)} \alpha + |q|^t \gamma e_{ij}^{(l)} \gamma^* \\
&= |q|^{- t} q a_+(l - \tfrac{1}{2}, -i + \tfrac{1}{2}, -j + \tfrac{1}{2}) \biggl(a_+(l - 1, i , j )\; e^{\left(l - 1\right)}_{i, j} 
+ a_-(l,i, j)\; e^{\left(l \right)}_{i, j} \biggr) \\
&\quad + |q|^{- t} q a_-(l + \tfrac{1}{2},-i + \tfrac{1}{2}, -j + \tfrac{1}{2}) \biggl(a_+(l, i , j )\; e^{\left(l\right)}_{i, j} 
+ a_-(l + 1,i, j)\; e^{\left(l + 1\right)}_{i, j} \biggr) \allowdisplaybreaks[2] \\
&\quad - |q|^t c_+ (l,-i,-j) \biggl(c_+ (l + \tfrac{1}{2}, i - \tfrac{1}{2} ,j + \tfrac{1}{2}) \, e_{i, j}^{\left(l + 1\right)} 
+ c_- (l + \tfrac{1}{2},i - \tfrac{1}{2},j + \tfrac{1}{2}) \, e_{i, j}^{\left(l \right)} \biggr) \\
&\quad - |q|^t c_-(l,-i,-j) \biggl(c_+ (l - \tfrac{1}{2}, i - \tfrac{1}{2} ,j + \tfrac{1}{2}) \, e_{i, j}^{\left(l\right)} 
+ c_- (l - \tfrac{1}{2},i - \tfrac{1}{2},j + \tfrac{1}{2}) \, e_{i, j}^{\left(l - 1\right)} \biggr).  
\end{align*} 
Similarly we obtain 
\begin{align*}
&\pi_t(\gamma)(e^{(l)}_{ij}) = |q|^t \gamma e^{(l)}_{ij} \alpha^* - |q|^{-t} q \alpha^* e^{(l)}_{ij} \gamma \\
&= |q|^{t} q^{-1} a_+(l,-i,-j) \biggl( c_+ (l + \tfrac{1}{2},i + \tfrac{1}{2},j + \tfrac{1}{2}) \, e_{i + 1, j}^{\left(l + 1 \right)} + 
c_- (l + \tfrac{1}{2},i + \tfrac{1}{2},j + \tfrac{1}{2}) \, e_{i + 1, j}^{\left(l\right)} \biggr) \\
&+ |q|^{t} q^{-1} a_-(l,-i,-j) \biggl( 
c_+ (l - \tfrac{1}{2},i + \tfrac{1}{2},j + \tfrac{1}{2}) \, e_{i + 1, j}^{\left(l \right)} + 
c_- (l - \tfrac{1}{2},i + \tfrac{1}{2},j + \tfrac{1}{2}) \, e_{i + 1, j}^{\left(l - 1 \right)} \biggr) \allowdisplaybreaks[2] \\
& + |q|^{-t} q c_+(l - \tfrac{1}{2}, -i - \tfrac{1}{2}, -j + \tfrac{1}{2}) \biggl(
a_+(l - 1, i + 1, j)\; e^{\left(l - 1\right)}_{i + 1, j} 
+ a_-(l,i + 1, j)\; e^{\left(l\right)}_{i + 1, j} \biggr) \\
&+ |q|^{-t} q c_-(l + \tfrac{1}{2}, -i - \tfrac{1}{2}, -j + \tfrac{1}{2})\; \biggl( 
a_+(l, i + 1, j)\; e^{\left(l\right)}_{i + 1, j} 
+ a_-(l + 1,i + 1, j)\; e^{\left(l + 1\right)}_{i + 1, j}\biggr) 
\end{align*}
and 
\begin{align*}
&\pi_t(\gamma^*)(e^{(l)}_{ij}) = |q|^{-t} \gamma^* e^{(l)}_{ij} \alpha - |q|^{t} q^{-1} \alpha e^{(l)}_{ij} \gamma^* \\
&= |q|^{-t} q a_+(l - \tfrac{1}{2}, -i + \tfrac{1}{2}, -j + \tfrac{1}{2}) \biggl(c_+(l - 1, i - 1, j)\; e^{\left(l - 1\right)}_{i - 1, j} 
+ c_-(l,i - 1, j)\; e^{\left(l \right)}_{i - 1, j} \biggr) \\
& + |q|^{-t} q a_-(l + \tfrac{1}{2},-i + \tfrac{1}{2}, -j + \tfrac{1}{2}) 
\biggl(c_+(l, i - 1, j)\; e^{\left(l\right)}_{i - 1, j} + c_-(l + 1,i - 1, j)\; e^{\left(l + 1 \right)}_{i - 1, j} \biggr) \allowdisplaybreaks[2]\\
& + |q|^{t} q^{-1} c_+ (l,-i,-j) \biggl(a_+ (l + \tfrac{1}{2}, i - \tfrac{1}{2} ,j + \tfrac{1}{2}) \, e_{i - 1, j}^{\left(l + 1\right)} + 
a_- (l + \tfrac{1}{2},i - \tfrac{1}{2},j + \tfrac{1}{2}) \, e_{i - 1, j}^{\left(l \right)} \biggr) \\
& + |q|^{t} q^{-1} c_-(l,-i,-j) \biggl(a_+ (l - \tfrac{1}{2}, i - \tfrac{1}{2} ,j + \tfrac{1}{2}) \, e_{i - 1, j}^{\left(l\right)} + 
a_- (l - \tfrac{1}{2},i - \tfrac{1}{2},j + \tfrac{1}{2}) \, e_{i - 1, j}^{\left(l - 1\right)} \biggr). 
\end{align*}
This may be written in the form 
\begin{align*}
\pi_t(\alpha)(e^{(l)}_{i,j}) &= a_1(t,l,i,j) e^{(l + 1)}_{i,j} + a_0(t,l,i,j) e^{(l)}_{i,j} + a_{-1}(t,l,i,j) e^{(l - 1)}_{i,j} \\
\pi_t(\alpha^*)(e^{(l)}_{i,j}) &= b_{1}(t,l,i,j) e^{(l + 1)}_{i,j} + b_0(t,l,i,j) e^{(l)}_{i,j} + b_{-1}(t,l,i,j) e^{(l - 1)}_{i,j} \\
\pi_t(\gamma)(e^{(l)}_{i,j}) &= c_1(t,l,i,j) e^{(l + 1)}_{i + 1,j} + c_0(t,l,i,j) e^{(l)}_{i + 1,j} + c_{-1}(t,l,i,j) e^{(l - 1)}_{i + 1,j} \\
\pi_t(\gamma^*)(e^{(l)}_{i,j}) &= d_1(t,l,i,j) e^{(l + 1)}_{i - 1,j} + d_0(t,l,i,j) e^{(l)}_{i - 1,j} + d_{-1}(t,l,i,j) e^{(l - 1)}_{i - 1,j} 
\end{align*}
where 
\begin{align*}
a_1(t,l,i,j) &= |q|^t 
q^{4l + 3} \, 
\frac{(1 - q^{2l + 2j + 2})^{\tfrac{1}{2}} (1 - q^{2l + 2i + 2})^{\tfrac{1}{2}} (1 - q^{2l - 2j + 2})^{\tfrac{1}{2}} 
(1 - q^{2l - 2i + 2})^{\tfrac{1}{2}}}
{(1 - q^{4l + 2})^{\tfrac{1}{2}} (1 - q^{4l + 4}) (1 - q^{4l + 6})^{\tfrac{1}{2}}} \\
&- |q|^{-t} q^{2l + 3} \, \frac{(1 - q^{2l - 2j + 2})^{\tfrac{1}{2}} 
(1 - q^{2l + 2i + 2})^{\tfrac{1}{2}}(1 - q^{2l + 2j + 2})^{\tfrac{1}{2}} (1 - q^{2l - 2i + 2})^{\tfrac{1}{2}}}
{(1 - q^{4l + 2})^{\tfrac{1}{2}} (1 - q^{4l + 4}) (1 - q^{4l + 6})^{\tfrac{1}{2}}} \allowdisplaybreaks[2] \\
a_0(t,l,i,j) &= |q|^t q^{2l - i - j} \, \frac{(1 - q^{2l + 2j + 2}) (1 - q^{2l + 2i + 2})}{ (1 - q^{4l + 2}) (1 - q^{4l + 4})} \\
&\qquad + |q|^t q^{2l + i + j} \frac{(1 - q^{2l - 2j}) (1 - q^{2l - 2i})}{(1 - q^{4l}) (1 - q^{4l+2})} \\
&\qquad + |q|^{-t} q^{2l + i - j + 2} \frac{(1 - q^{2l + 2j})(1 - q^{2l - 2i})}{(1 - q^{4l}) (1 - q^{4l + 2})} \\
&\qquad + |q|^{-t} q^{2l - i + j + 2} \, \frac{(1 - q^{2l - 2j + 2}) (1 - q^{2l + 2i + 2})}{(1 - q^{4l + 2}) (1 - q^{4l + 4})} \allowdisplaybreaks[2]\\
a_{-1}(t,l,i,j) &= |q|^t q^{-1} 
\frac{(1 - q^{2l - 2j})^{\tfrac{1}{2}} (1 - q^{2l - 2i})^{\tfrac{1}{2}}(1 - q^{2l + 2j})^{\tfrac{1}{2}} 
(1 - q^{2l + 2i})^{\tfrac{1}{2}}}{(1 - q^{4l})(1 - q^{4l+2})^{\tfrac{1}{2}} (1 - q^{4l - 2})^{\tfrac{1}{2}}} \\
& - |q|^{-t} q^{2l + 1} \, \frac{(1 - q^{2l + 2j})^{\tfrac{1}{2}} 
(1 - q^{2l - 2i})^{\tfrac{1}{2}}(1 - q^{2l - 2j})^{\tfrac{1}{2}} (1 - q^{2l + 2i})^{\tfrac{1}{2}}}{(1 - q^{4l}) (1 - q^{4l + 2})^{\tfrac{1}{2}} (1 - q^{4l - 2})^{\tfrac{1}{2}}} 
\end{align*}
and
\begin{align*}
b_1(t,l,i,j) &= |q|^{-t} q\frac{(1 - q^{2l - 2j + 2})^{\tfrac{1}{2}} (1 - q^{2l - 2i + 2})^{\tfrac{1}{2}}(1 - q^{2l + 2j + 2})^{\tfrac{1}{2}} (1 - q^{2l + 2i + 2})^{\tfrac{1}{2}}}
{(1 - q^{4l + 2})^{\tfrac{1}{2}} (1 - q^{4l+4}) (1 - q^{4l+6})^{\tfrac{1}{2}}} \\
& - |q|^t q^{2l + 1} \, \frac{(1 - q^{2l + 2j + 2})^{\tfrac{1}{2}} (1 - q^{2l - 2i + 2})^{\tfrac{1}{2}}(1 - q^{2l - 2j + 2})^{\tfrac{1}{2}} 
(1 - q^{2l + 2i + 2})^{\tfrac{1}{2}}}{(1 - q^{4l + 2})^{\tfrac{1}{2}} (1 - q^{4l + 4}) (1 - q^{4l + 6})^{\tfrac{1}{2}}} \allowdisplaybreaks[2]\\
b_0(t,l,i,j) &= |q|^{-t} q^{2l - i - j + 2} \, \frac{(1 - q^{2l + 2j}) (1 - q^{2l + 2i})}{ (1 - q^{4l}) (1 - q^{4l + 2})} \\
&\qquad + |q|^{-t} q^{2l + i + j + 2} \frac{(1 - q^{2l - 2j + 2}) (1 - q^{2l - 2i + 2})}{(1 - q^{4l + 2}) (1 - q^{4l+4})} \\
&\qquad + |q|^t q^{2l + i - j} \, \frac{(1 - q^{2l + 2j + 2})(1 - q^{2l - 2i + 2})}{(1 - q^{4l + 2}) (1 - q^{4l + 4})} \\
&\qquad + |q|^t q^{2l - i + j} \, \frac{(1 - q^{2l - 2j})(1 - q^{2l + 2i})}{(1 - q^{4l}) (1 - q^{4l + 2})} \allowdisplaybreaks[2] \\
b_{-1}(t,l,i,j) &= |q|^{-t} q^{4l + 1} \, 
\frac{(1 - q^{2l + 2j})^{\tfrac{1}{2}} (1 - q^{2l + 2i})^{\tfrac{1}{2}}(1 - q^{2l - 2j})^{\tfrac{1}{2}} (1 - q^{2l - 2i})^{\tfrac{1}{2}}}
{(1 - q^{4l}) (1 - q^{4l + 2})^{\tfrac{1}{2}} (1 - q^{4l - 2})^{\tfrac{1}{2}}} \\
&- |q|^t q^{2l - 1} \, \frac{(1 - q^{2l - 2j})^{\tfrac{1}{2}} (1 - q^{2l + 2i})^{\tfrac{1}{2}}
(1 - q^{2l + 2j})^{\tfrac{1}{2}} (1 - q^{2l - 2i})^{\tfrac{1}{2}}}
{(1 - q^{4l})(1 - q^{4l + 2})^{\tfrac{1}{2}} (1 - q^{4l - 2})^{\tfrac{1}{2}}},
\end{align*} 
similarly
\begin{align*}
c_1(t,l,i,j) &= - |q|^t q^{3l - i + 1} \, \frac{(1 - q^{2l + 2j + 2})^{\tfrac{1}{2}} (1 - q^{2l + 2i + 2})^{\tfrac{1}{2}} 
(1 - q^{2l - 2j + 2})^{\tfrac{1}{2}} (1 - q^{2l + 2i + 4})^{\tfrac{1}{2}}}
{(1 - q^{4l + 2})^{\tfrac{1}{2}} (1 - q^{4l + 4}) (1 - q^{4l + 6})^{\tfrac{1}{2}}} \\
&+ |q|^{-t} q^{l - i + 1} \, \frac{(1 - q^{2l - 2j + 2})^{\tfrac{1}{2}} 
(1 - q^{2l + 2i + 2})^{\tfrac{1}{2}} (1 - q^{2l + 2j + 2})^{\tfrac{1}{2}} (1 - q^{2l + 2i + 4})^{\tfrac{1}{2}}}
{(1 - q^{4l + 2})^{\tfrac{1}{2}} (1 - q^{4l + 4}) (1 - q^{4l+6})^{\tfrac{1}{2}}}  \allowdisplaybreaks[2]\\
c_0(t,l,i,j) &= |q|^t q^{3l - j + 1} \, \frac{(1 - q^{2l + 2j + 2})(1 - q^{2l + 2i + 2})^{\tfrac{1}{2}}(1 - q^{2l - 2i})^{\tfrac{1}{2}}}
{ (1 - q^{4l + 2})(1 - q^{4l + 4})} \\
&\qquad - |q|^t q^{l + j - 1} \frac{(1 - q^{2l - 2j}) (1 - q^{2l - 2i})^{\tfrac{1}{2}} (1 - q^{2l + 2i + 2})^{\tfrac{1}{2}}}{(1 - q^{4l}) (1 - q^{4l+2})} \\
&\qquad - |q|^{-t} q^{l - j + 1} \, \frac{(1 - q^{2l + 2j}) (1 - q^{2l - 2i})^{\tfrac{1}{2}}(1 - q^{2l + 2i + 2})^{\tfrac{1}{2}}}
{(1 - q^{4l})(1 - q^{4l + 2})} \\
&\qquad + |q|^{-t} q^{3l + j + 3} \, \frac{(1 - q^{2l - 2j + 2}) (1 - q^{2l + 2i + 2})^{\tfrac{1}{2}}(1 - q^{2l - 2i})^{\tfrac{1}{2}}}
{(1 - q^{4l + 2}) (1 - q^{4l + 4})} \allowdisplaybreaks[2]\\
c_{-1}(t,l,i,j) &= |q|^t q^{l + i - 1} 
\frac{(1 - q^{2l - 2j})^{\tfrac{1}{2}} (1 - q^{2l - 2i})^{\tfrac{1}{2}}(1 - q^{2l + 2j})^{\tfrac{1}{2}} (1 - q^{2l - 2i - 2})^{\tfrac{1}{2}}}
{(1 - q^{4l}) (1 - q^{4l + 2})^{\tfrac{1}{2}} (1 - q^{4l - 2})^{\tfrac{1}{2}}} \\
& - |q|^{-t} q^{3l + i + 1} \, \frac{(1 - q^{2l + 2j})^{\tfrac{1}{2}} (1 - q^{2l - 2i})^{\tfrac{1}{2}}
(1 - q^{2l - 2j})^{\tfrac{1}{2}} (1 - q^{2l - 2i - 2})^{\tfrac{1}{2}}}
{(1 - q^{4l})(1 - q^{4l + 2})^{\tfrac{1}{2}} (1 - q^{4l - 2})^{\tfrac{1}{2}}}, 
\end{align*} 
and
\begin{align*}
d_1(t,l,i,j) &= |q|^{-t} q^{l + i + 1} 
\frac{(1 - q^{2l - 2j + 2})^{\tfrac{1}{2}} (1 - q^{2l - 2i + 2})^{\tfrac{1}{2}} (1 - q^{2l + 2j + 2})^{\tfrac{1}{2}} 
(1 - q^{2l - 2i + 4})^{\tfrac{1}{2}}}
{(1 - q^{4l + 2})^{\tfrac{1}{2}} (1 - q^{4l+4}) (1 - q^{4l + 6})^{\tfrac{1}{2}}} \\
&- |q|^t q^{3l + i + 1} \, \frac{(1 - q^{2l + 2j + 2})^{\tfrac{1}{2}} (1 - q^{2l - 2i + 2})^{\tfrac{1}{2}}(1 - q^{2l - 2j + 2})^{\tfrac{1}{2}} (1 - q^{2l - 2i + 4})^{\tfrac{1}{2}}}
{(1 - q^{4l + 2})^{\tfrac{1}{2}} (1 - q^{4l + 4}) (1 - q^{4l + 6})^{\tfrac{1}{2}}} \allowdisplaybreaks[2] \\
d_0(t,l,i,j) &= |q|^{- t} q^{3l - j + 1} \, 
\frac{(1 - q^{2l + 2j}) (1 - q^{2l + 2i})^{\tfrac{1}{2}}(1 - q^{2l - 2i + 2})^{\tfrac{1}{2}}}{ (1 - q^{4l}) (1 - q^{4l + 2})} \\
&\qquad- |q|^{- t} q^{l + j + 1} \frac{(1 - q^{2l - 2j + 2}) (1 - q^{2l - 2i + 2})^{\tfrac{1}{2}}(1 - q^{2l + 2i})^{\tfrac{1}{2}}}
{(1 - q^{4l + 2}) (1 - q^{4l + 4})} \\
&\qquad - |q|^t q^{l - j - 1} \, \frac{(1 - q^{2l + 2j + 2})(1 - q^{2l - 2i + 2})^{\tfrac{1}{2}}(1 - q^{2l + 2i})^{\tfrac{1}{2}}}
{(1 - q^{4l + 2}) (1 - q^{4l + 4})} \\
&\qquad + |q|^t q^{3l + j - 1} \, \frac{(1 - q^{2l - 2j})(1 - q^{2l + 2i})^{\tfrac{1}{2}} (1 - q^{2l - 2i + 2})^{\tfrac{1}{2}}}
{(1 - q^{4l}) (1 - q^{4l + 2})} \allowdisplaybreaks[2] \\
d_{-1}(t,l,i,j) &= - |q|^{- t} q^{3l - i + 1} 
\, \frac{(1 - q^{2l + 2j})^{\tfrac{1}{2}} (1 - q^{2l + 2i})^{\tfrac{1}{2}}(1 - q^{2l - 2j})^{\tfrac{1}{2}} (1 - q^{2l + 2i - 2})^{\tfrac{1}{2}}}
{(1 - q^{4l}) (1 - q^{4l + 2})^{\tfrac{1}{2}} (1 - q^{4l - 2})^{\tfrac{1}{2}}} \\
& + |q|^t q^{l - i - 1} \, \frac{(1 - q^{2l - 2j})^{\tfrac{1}{2}} (1 - q^{2l + 2i})^{\tfrac{1}{2}}(1 - q^{2l + 2j})^{\tfrac{1}{2}} 
(1 - q^{2l + 2i - 2})^{\tfrac{1}{2}}}
{(1 - q^{4l}) (1 - q^{4l + 2})^{\tfrac{1}{2}} (1 - q^{4l - 2})^{\tfrac{1}{2}}}.
\end{align*} 
Let us now set 
$$
m(t,l) = \frac{|q|^{-t} q^2 - |q|^t q^{2l}}{|q|^t - |q|^{-t} q^{2l + 2}} = \frac{q^2 - |q|^{2t} q^{2l}}{|q|^{2t} - q^{2l + 2}}
$$
for $ t \in [0,1] $ and $ l \in \mathbb{N} $. Note that $ m(1,l) = 1 $ for all $ l $ if we interpret $ m(1,0) = 1 $. 
We define 
\begin{align*}
A_1(t,l,i) &= m(t, l + 1)^{-\frac{1}{2}}\, a_1(t,l,i,0) \\
A_0(t,l,i) &= a_0(t,l,i,0) \\
A_{-1}(t,l,i) &= m(t, l)^{\frac{1}{2}}\, a_{-1}(t,l,i,0), 
\end{align*}
and by rescaling $ b_k(t,l,i,0), c_k(t,l,i,0) $ and $ d_k(t,l,i,0) $ for $ k = -1,0,1 $
in the same way we obtain constants $ B_k(t,l,i), C_k(t,l,i) $ and $ D_k(t,l,i) $. 
Inspection of the formulas above shows that the expressions $ X_1(t,0,0) $ for $ X = A, B, C, D $ are
well-defined and depend continuously on $ t \in [0,1] $ although $ m(0, 1) = 0 $. 
\begin{lemma} \label{lemma1}
Let $ q \in (-1,1) \setminus \{0\} $. For $ t \in [0,1] $ the formulas 
\begin{align*}
\omega_t(\alpha)(e^{(l)}_{i,0}) &= A_1(t,l,i) e^{(l + 1)}_{i,0} + A_0(t,l,i) e^{(l)}_{i,0} + A_{-1}(t,l,i) e^{(l - 1)}_{i,0} \\
\omega_t(\alpha^*)(e^{(l)}_{i,0}) &= B_{1}(t,l,i) e^{(l + 1)}_{i,0} + B_0(t,l,i) e^{(l)}_{i,0} + B_{-1}(t,l,i) e^{(l - 1)}_{i,0} \\
\omega_t(\gamma)(e^{(l)}_{i,0}) &= C_1(t,l,i) e^{(l + 1)}_{i + 1,0} + C_0(t,l,i) e^{(l)}_{i + 1,0} + C_{-1}(t,l,i) e^{(l - 1)}_{i + 1,0} \\
\omega_t(\gamma^*)(e^{(l)}_{i,0}) &= D_1(t,l,i) e^{(l + 1)}_{i - 1,0} + D_0(t,l,i) e^{(l)}_{i - 1,0} + D_{-1}(t,l,i) e^{(l - 1)}_{i - 1,0} 
\end{align*}
define a $ * $-homomorphism $ \omega_t: C(G_q) \rightarrow \LH(L^2(\E_0)) $. 
\end{lemma}
\proof The main point is to show that $ \omega_t $ is compatible with the $ * $-structures. 
In order to prove $ \omega_t(\alpha)^* = \omega_t(\alpha^*) $ we have to verify 
\begin{align*}
A_1(t,l,i) &= B_{-1}(t,l + 1,i) \\
A_0(t,l,i) &= B_0(t,l,i) \\
A_{-1}(t,l,i) &= B_1(t,l - 1,i).
\end{align*}
We obtain
\begin{align*}
m(t, l + 1)^{\frac{1}{2}} A_1(t,l,i) &= (|q|^t q^{4l + 3} - |q|^{-t} q^{2l + 3}) \times \\
&\qquad\qquad \frac{(1 - q^{2l + 2})(1 - q^{2l + 2i + 2})^{\tfrac{1}{2}} (1 - q^{2l - 2i + 2})^{\tfrac{1}{2}}}
{(1 - q^{4l + 2})^{\tfrac{1}{2}} (1 - q^{4l + 4}) (1 - q^{4l + 6})^{\tfrac{1}{2}}} 
\end{align*}
and 
\begin{align*}
m(t,l + 1)^{\frac{1}{2}} &B_{-1}(t,l + 1,i,0) = \frac{|q|^t q^{4l + 3} - |q|^{-t} q^{2l + 3}}{|q|^{-t} q^{4l + 5} - |q|^t q^{2l + 1}} \times \\
&(|q|^{-t} q^{4l + 5} - |q|^t q^{2l + 1})
\frac{(1 - q^{2l + 2})(1 - q^{2l + 2i + 2})^{\tfrac{1}{2}}(1 - q^{2l - 2i + 2})^{\tfrac{1}{2}}}
{(1 - q^{4l + 4}) (1 - q^{4l + 6})^{\tfrac{1}{2}} (1 - q^{4l + 2})^{\tfrac{1}{2}}}.
\end{align*}
Similarly we find
\begin{align*}
A_0(t,l,i) &= (|q|^t q^{2l - i} + |q|^{-t} q^{2l - i + 2})
\, \frac{(1 - q^{2l + 2}) (1 - q^{2l + 2i + 2})}{ (1 - q^{4l + 2}) (1 - q^{4l + 4})} \\
&\qquad+ (|q|^t q^{2l + i} + |q|^{-t} q^{2l + i + 2}) \frac{(1 - q^{2l}) (1 - q^{2l - 2i})}{(1 - q^{4l}) (1 - q^{4l+2})} \\
&= (|q|^t q^{2l - i} + |q|^{-t} q^{2l - i + 2})
\, \frac{(1 - q^{2l + 2i + 2})}{ (1 + q^{2l + 2}) (1 - q^{4l + 2})} \\
&\qquad+ (|q|^t q^{2l + i} + |q|^{-t} q^{2l + i + 2}) \frac{(1 - q^{2l - 2i})}{(1 + q^{2l}) (1 - q^{4l+2})} 
\end{align*}
and hence 
\begin{align*}
(1 + q^{2l})&(1 + q^{2l + 2})(1 - q^{4l+2}) A_0(t,l,i) \\
&= (|q|^t q^{2l - i} + |q|^{-t} q^{2l - i + 2})(1 - q^{2l + 2i + 2})(1 + q^{2l}) \\
&\qquad+ (|q|^t q^{2l + i} + |q|^{-t} q^{2l + i + 2}) (1 - q^{2l - 2i})(1 + q^{2l + 2}) \\
&= (|q|^t q^{2l - i} + |q|^{-t} q^{2l - i + 2})(1 - q^{2l + 2i + 2} + q^{2l} - q^{4l +2i + 2}) \\
&\qquad+ (|q|^t q^{2l + i} + |q|^{-t} q^{2l + i + 2}) (1 - q^{2l - 2i} + q^{2l + 2} - q^{4l - 2i + 2}) \allowdisplaybreaks[2] \\
&= |q|^t (q^{2l - i} - q^{4l + i + 2} + q^{4l - i} - q^{6l + i + 2} \\
&\qquad + q^{2l + i} - q^{4l - i} + q^{4l + i + 2} - q^{6l - i + 2}) \\
&\quad + |q|^{-t} (q^{2l - i + 2} - q^{4l + i + 4} + q^{4l - i + 2} - q^{6l +i + 4} \\
&\qquad+ q^{2l + i + 2} - q^{4l - i + 2} + q^{4l + i + 4} - q^{6l - i + 4}) \\
&= |q|^t (q^{2l - i} - q^{6l + i + 2} + q^{2l + i} - q^{6l - i + 2}) \\
&\quad + |q|^{-t} (q^{2l - i + 2} - q^{6l +i + 4} + q^{2l + i + 2} - q^{6l - i + 4}).
\end{align*}
Conversely, 
\begin{align*}
B_0(t,l,i) &= (|q|^{-t} q^{2l + i + 2} + |q|^t q^{2l + i}) 
\frac{(1 - q^{2l + 2}) (1 - q^{2l - 2i + 2})}{(1 - q^{4l + 2}) (1 - q^{4l+4})} \\
&\qquad + (|q|^{-t} q^{2l - i + 2} + |q|^t q^{2l - i}) \, \frac{(1 - q^{2l}) (1 - q^{2l + 2i})}{ (1 - q^{4l}) (1 - q^{4l + 2})} \\
&= (|q|^{-t} q^{2l + i + 2} + |q|^t q^{2l + i}) 
\frac{(1 - q^{2l - 2i + 2})}{(1 - q^{4l + 2}) (1 + q^{2l + 2})} \\
&\qquad + (|q|^{-t} q^{2l - i + 2} + |q|^t q^{2l - i}) \, \frac{(1 - q^{2l + 2i})}{ (1 + q^{2l}) (1 - q^{4l + 2})}
\end{align*}
and hence 
\begin{align*}
(1 + q^{2l})&(1 + q^{2l + 2})(1 - q^{4l+2}) B_0(t,l,i) \\
&= (|q|^{-t} q^{2l + i + 2} + |q|^t q^{2l + i}) 
(1 - q^{2l - 2i + 2}) (1 + q^{2l}) \\
&\qquad + (|q|^{-t} q^{2l - i + 2} + |q|^t q^{2l - i}) (1 - q^{2l + 2i}) (1 + q^{2l + 2}) \\
&= (|q|^{-t} q^{2l + i + 2} + |q|^t q^{2l + i}) 
(1 - q^{2l - 2i + 2} + q^{2l} - q^{4l - 2i + 2}) \\
&\qquad + (|q|^{-t} q^{2l - i + 2} + |q|^t q^{2l - i}) (1 - q^{2l + 2i} + q^{2l + 2} - q^{4l + 2i + 2}) \allowdisplaybreaks[2] \\
&= |q|^t (q^{2l + i} - q^{4l - i + 2} + q^{4l + i} - q^{6l - i + 2} \\
&\qquad + q^{2l - i} - q^{4l + i} + q^{4l - i + 2} - q^{6l + i + 2}) \\
&\quad + |q|^{-t} (q^{2l + i + 2} - q^{4l - i + 4} + q^{4l + i + 2} - q^{6l - i + 4} \\
&\qquad + q^{2l - i + 2} - q^{4l + i + 2} + q^{4l - i + 4} - q^{6l + i + 4}) \\
&= |q|^t (q^{2l + i} - q^{6l - i + 2} + q^{2l - i} - q^{6l + i + 2}) \\
&\qquad + |q|^{-t} (q^{2l + i + 2} - q^{6l - i + 4} + q^{2l - i + 2} - q^{6l + i + 4}).
\end{align*}
Finally,
\begin{align*}
m(t,l)^{\frac{1}{2}} &A_{-1}(t,l,i,0) = \frac{|q|^t q^{2l - 1} - |q|^{-t} q}{|q|^{-t} q^{2l + 1} - |q|^t q^{-1}} \times \\
&(|q|^t q^{-1} - |q|^{-t} q^{2l + 1})
\frac{(1 - q^{2l}) (1 - q^{2l - 2i})^{\tfrac{1}{2}} (1 - q^{2l + 2i})^{\tfrac{1}{2}}}
{(1 - q^{4l})(1 - q^{4l+2})^{\tfrac{1}{2}} (1 - q^{4l - 2})^{\tfrac{1}{2}}} 
\end{align*}
and 
\begin{align*}
m(t,l)^{\frac{1}{2}} B_1(t,l - 1,i,0) &= (|q|^{-t} q - |q|^t q^{2l - 1}) \frac{(1 - q^{2l}) (1 - q^{2l - 2i})^{\tfrac{1}{2}} (1 - q^{2l + 2i})^{\tfrac{1}{2}}}
{(1 - q^{4l - 2})^{\tfrac{1}{2}} (1 - q^{4l}) (1 - q^{4l + 2})^{\tfrac{1}{2}}}.
\end{align*}
One verifies $ \omega_t(\gamma)^* = \omega_t(\gamma^*) $ in a similar fashion. 
It is then easy to check that the operators $ \omega_t(\alpha) $ and $ \omega_t(\gamma) $ 
satisfy the defining relations of $ C(G_q) $. \qed 
\begin{lemma} \label{lemma2} 
Let $ q \in (-1,1) \setminus \{0\} $. For $ l \rightarrow \infty $ the expressions 
\begin{align*}
|a_1(1, l, i, \pm 1) - A_1(t, l, i)|, |a_0(1, l, i, \pm 1)|, |A_0(t, l, i)|, |a_{-1}(1, l, i, \pm 1) - A_{-1}(t, l, i)| 
\end{align*}
and 
\begin{align*}
|c_1(1, l, i, \pm 1) - C_1(t, l, i)|, |c_0(1, l, i, \pm 1)|, |C_0(t, l, i)|, |c_{-1}(1, l, i, \pm 1) - C_{-1}(t, l, i)| 
\end{align*}
tend to zero uniformly for $ t \in [0,1] $ and independently of $ i $. 
\end{lemma} 
\proof Since the constants $ a_k(1,l,i,j) $ and $ c_k(1,l,i,j) $ for $ k = -1,0,1 $ are symmetric 
in the variable $ j $ we may restrict attention to the case $ j = 1 $. \\
The estimates involving $ A_1, A_0 $ and $ a_0 $ are easy. For $ A_{-1} $ it suffices to consider
\begin{align*}
\biggl|m(t,l)^{\frac{1}{2}} |q|^t q^{-1} 
\frac{(1 - q^{2l})(1 - q^{2l - 2i})^{\tfrac{1}{2}}(1 - q^{2l + 2i})^{\tfrac{1}{2}}}
{(1 - q^{4l})(1 - q^{4l+2})^{\tfrac{1}{2}} (1 - q^{4l - 2})^{\tfrac{1}{2}}} \\
- |q| q^{-1} 
\frac{(1 - q^{2l - 2})^{\tfrac{1}{2}} (1 - q^{2l - 2i})^{\tfrac{1}{2}}(1 - q^{2l + 2})^{\tfrac{1}{2}} 
(1 - q^{2l + 2i})^{\tfrac{1}{2}}}{(1 - q^{4l})(1 - q^{4l+2})^{\tfrac{1}{2}} (1 - q^{4l - 2})^{\tfrac{1}{2}}} \biggr|
\end{align*}
which reduces to 
\begin{align*}
\biggl|\frac{(|q|^t q^{2l} - |q|^{-t} q^2)^{\frac{1}{2}}}{(|q|^{-t} q^{2l + 2} - |q|^t)^{\frac{1}{2}}} |q|^t 
(1 - q^{2l})- |q| (1 - q^{2l - 2})^{\tfrac{1}{2}} (1 - q^{2l + 2})^{\tfrac{1}{2}} \biggr|.
\end{align*}
It is enough to estimate
\begin{align*}
&\biggl|\frac{|q|^t q^{2l} - |q|^{-t} q^2}{|q|^{-t} q^{2l + 2} - |q|^t} |q|^{2t} 
(1 - q^{2l})^2 - |q|^2 (1 - q^{2l - 2})(1 - q^{2l + 2}) \biggr| \\
&= \biggl| \frac{|q|^{2t} q^{2l} - q^2}{|q|^{-2t} q^{2l + 2} - 1} 
(1 - q^{2l})^2 - |q|^2 (1 - q^{2l - 2})(1 - q^{2l + 2}) \biggr|. 
\end{align*}
We may estimate this expression by 
\begin{align*}
&\biggl| \frac{|q|^{2t} q^{2l}}{|q|^{-2t} q^{2l + 2} - 1} 
(1 - q^{2l})^2 \biggr| \\
&\qquad + \biggl|\frac{q^2}{1 - |q|^{-2t} q^{2l + 2}} 
(1 - q^{2l})^2 - |q|^2 (1 - q^{2l - 2})(1 - q^{2l + 2}) \biggr|,  
\end{align*}
and both terms converge to zero for $ l \rightarrow \infty $. The 
remaining assertions are verified in a similar fashion. \qed \\
In the sequel we write $ \sgn(q) $ for the sign of $ q $, that is, $ \sgn(q) = 1 $ if $ q > 0 $ and 
$ \sgn(q) = -1 $ if $ q < 0 $. 
\begin{lemma} \label{lemma3}
Let $ q \in (-1,1) \setminus \{0\} $. We have 
\begin{align*}
A_1(0,l,i) &= a_1(1,l,i, \pm 1) \\ 
A_0(0,l,i) &= \sgn(q)\, a_0(1,l,i, \pm 1) \\
A_{-1}(0,l,i) &= a_{-1}(1,l,i,\pm 1)
\end{align*}
and similarly
\begin{align*}
C_1(0,l,i) &= c_1(1,l,i, \pm 1) \\ 
C_0(0,l,i) &= \sgn(q)\, c_0(1,l,i, \pm 1) \\
C_{-1}(0,l,i) &= c_{-1}(1,l,i, \pm 1)
\end{align*}
for $ l > 0 $ and all $ i $. 
\end{lemma}
\proof Since the coefficients $ a_k(1,l,i,j) $ and $ c_k(1,l,i,j) $ for $ k = -1,0,1 $ are symmetric 
in the variable $ j $ it suffices again to consider the case $ j = 1 $. We have 
\begin{align*}
A_1(0,l,i) &= m(0,l + 1)^{-\frac{1}{2}}(q^{4l + 3} - q^{2l + 3}) \times \\
&\qquad \frac{(1 - q^{2l + 2})
(1 - q^{2l + 2i + 2})^{\tfrac{1}{2}}(1 - q^{2l - 2i + 2})^{\tfrac{1}{2}}}
{(1 - q^{4l + 2})^{\tfrac{1}{2}} (1 - q^{4l + 4}) (1 - q^{4l + 6})^{\tfrac{1}{2}}} 
\end{align*}
and 
\begin{align*}
a_1(1,l,i,1) &= (|q| q^{4l + 3} - |q|^{-1} q^{2l + 3}) \times \\
&\qquad \, \frac{(1 - q^{2l})^{\tfrac{1}{2}} 
(1 - q^{2l + 2i + 2})^{\tfrac{1}{2}}(1 - q^{2l + 4})^{\tfrac{1}{2}} (1 - q^{2l - 2i + 2})^{\tfrac{1}{2}}}
{(1 - q^{4l + 2})^{\tfrac{1}{2}} (1 - q^{4l + 4}) (1 - q^{4l + 6})^{\tfrac{1}{2}}}.
\end{align*}
Then 
\begin{align*}
m(0,l + 1)^{-\frac{1}{2}}&(q^{2l + 1} - q) (1 - q^{2l + 2}) = 
\frac{(1 - q^{2l + 4})^{\frac{1}{2}}}{(q^2 - q^{2l + 2})^{\frac{1}{2}}} (q^{2l + 1} - q) (1 - q^{2l + 2}) \\
&= \frac{(1 - q^{2l + 4})^{\frac{1}{2}}}{(1 - q^{2l})^{\frac{1}{2}}} |q|^{-1} (-q) (1 - q^{2l}) (1 - q^{2l + 2}) \\
&= (1 - q^{2l + 4})^{\frac{1}{2}} |q|^{-1} (1 - q^{2l})^{\frac{1}{2}} (q^{2l + 3} - q) \\
&= (|q| q^{2l + 1} - |q|^{-1} q) (1 - q^{2l})^{\tfrac{1}{2}} (1 - q^{2l + 4})^{\tfrac{1}{2}}
\end{align*}
yields the first claim. Next we note 
\begin{align*}
A_0(0,l,i) &= (q^{2l - i} + q^{2l - i + 2}) \frac{(1 - q^{2l + 2}) (1 - q^{2l + 2i + 2})}{ (1 - q^{4l + 2}) (1 - q^{4l + 4})} \\
&\qquad + (q^{2l + i} + q^{2l + i + 2}) \frac{(1 - q^{2l}) (1 - q^{2l - 2i})}{(1 - q^{4l}) (1 - q^{4l+2})} 
\end{align*}
and 
\begin{align*}
a_0(1,l,i,1) &= (|q| q^{2l - i - 1}(1 - q^{2l + 4}) + |q|^{-1} q^{2l - i + 3} (1 - q^{2l})) 
\frac{1 - q^{2l + 2i + 2}}{ (1 - q^{4l + 2}) (1 - q^{4l + 4})} \\
&+ (|q| q^{2l + i + 1} (1 - q^{2l - 2}) + |q|^{-1} q^{2l + i + 1} (1 - q^{2l + 2})) \frac{1 - q^{2l - 2i}}{(1 - q^{4l}) (1 - q^{4l + 2})}. 
\end{align*}
Since 
\begin{align*}
\sgn(q) (q^{2l - i} + q^{2l - i + 2}) &(1 - q^{2l + 2}) = \sgn(q) (q^{2l - i} + q^{2l - i + 2} - q^{4l - i + 2} - q^{4l - i + 4}) \\
&= |q|q^{-1} (q^{2l - i} - q^{4l -i + 4}) + |q|^{-1}q (q^{2l - i + 2} - q^{4l - i + 2}) \\
&= |q| q^{2l - i - 1}(1 - q^{2l + 4}) + |q|^{-1} q^{2l - i + 3} (1 - q^{2l}) 
\end{align*}
and 
\begin{align*}
\sgn(q) (q^{2l + i} + q^{2l + i + 2})&(1 - q^{2l}) = \sgn(q) (q^{2l + i} + q^{2l + i + 2} - q^{4l + i} - q^{4l + i + 2}) \\
&= |q| q^{-1} (q^{2l + i + 2} - q^{4l + i}) + |q|^{-1} q (q^{2l + i} - q^{4l + i + 2}) \\
&= |q| q^{2l + i + 1} (1 - q^{2l - 2}) + |q|^{-1} q^{2l + i + 1} (1 - q^{2l + 2})
\end{align*}
we obtain the second assertion. 
Moreover, 
\begin{align*}
A_{-1}(0,l,i) &= m(0,l)^{\frac{1}{2}} (q^{-1} - q^{2l + 1}) (1 - q^{2l}) 
\frac{(1 - q^{2l - 2i})^{\tfrac{1}{2}} (1 - q^{2l + 2i})^{\tfrac{1}{2}}}{(1 - q^{4l}) 
(1 - q^{4l + 2})^{\tfrac{1}{2}} (1 - q^{4l - 2})^{\tfrac{1}{2}}} 
\end{align*}
and 
\begin{align*}
a_{-1}(1,l,i,1) &= (|q| q^{-1} - |q|^{-1} q^{2l + 1}) \times \\
&\qquad \frac{(1 - q^{2l - 2})^{\tfrac{1}{2}} (1 - q^{2l - 2i})^{\tfrac{1}{2}}(1 - q^{2l + 2})^{\tfrac{1}{2}} 
(1 - q^{2l + 2i})^{\tfrac{1}{2}}}{(1 - q^{4l})(1 - q^{4l+2})^{\tfrac{1}{2}} (1 - q^{4l - 2})^{\tfrac{1}{2}}}.  
\end{align*}
Hence 
\begin{align*}
m(0,l)^{\frac{1}{2}}(q^{-1} - q^{2l + 1})(1 - q^{2l}) &= \frac{(q^2 - q^{2l})^{\tfrac{1}{2}}}{(1 - q^{2l + 2})^{\tfrac{1}{2}}} 
(q^{-1} - q^{2l + 1})(1 - q^{2l}) \\
&= \frac{(1 - q^{2l - 2})^{\tfrac{1}{2}}}{(1 - q^{2l + 2})^{\tfrac{1}{2}}} |q|q^{-1}
(1 - q^{2l + 2})(1 - q^{2l}) \\
&= (1 - q^{2l - 2})^{\tfrac{1}{2}} |q|q^{-1}
(1 - q^{2l + 2})^{\tfrac{1}{2}} (1 - q^{2l}) \\
&= (|q| q^{-1} - |q|^{-1} q^{2l + 1}) (1 - q^{2l - 2})^{\tfrac{1}{2}} (1 - q^{2l + 2})^{\tfrac{1}{2}} 
\end{align*}
yields the third claim. The remaining assertions are verified in the same way. \qed \\
We need some further constructions. Recall that $ C(\E_k) $ for $ k \in \mathbb{Z} $ is a $ \DD(G_q) $-equivariant Hilbert $ C(G_q/T) $-module
in a natural way. Left multiplication yields a $ \DD(G_q) $-equivariant $ * $-homomorphism $ \psi: C(G_q/T) \rightarrow \KH(C(\E_k)) $, 
and $ (C(\E_k), \psi, 0) $ defines a class $ [[\E_k]] $ in
$ KK^{\DD(G_q)}(C(G_q/T), C(G_q/T)) $. Moreover $ [[\E_m]] \circ [[\E_n]] = [[\E_{m + n}]] $ for all $ m,n \in \mathbb{Z} $. \\
For $ k \in \mathbb{Z} $ we define $ [D_k] \in KK^{\DD(G_q)}(C(G_q/T), \mathbb{C}) $ by
$$
[D_k] = [[\E_k]] \circ [D]
$$
where $ [D] \in KK^{\DD(G_q)}(C(G_q/T), \mathbb{C}) $ is the element obtained in 
proposition \ref{DiracYD}. Remark that $ [D_0] = [D] $ since $ [[\E_0]] = 1 $. \\
Evidently, the unit homomorphism $ u: \mathbb{C} \rightarrow C(G_q/T) $ induces a class $ [u] $ in $ KK^{\DD(G_q)}(\mathbb{C}, C(G_q/T)) $.
We define $ [\E_k] $ in $ KK^{\DD(G_q)}(\mathbb{C}, C(G_q/T)) $ by restricting $ [[\E_k]] $ along $ u $,
or equivalently, by taking the product
$$
[\E_k] = [u] \circ [[\E_k]].
$$
In the sequel we will interested in the elements 
$ \alpha_q \in KK^{\DD(G_q)}(C(G_q/T), \mathbb{C} \oplus \mathbb{C}) $ and $ \beta_q 
\in KK^{\DD(G_q)}(\mathbb{C} \oplus \mathbb{C}, C(G_q/T)) $ given by 
$$
\alpha_q = [D_0] \oplus [D_{-1}], \qquad \beta_q = (-[\E_1]) \oplus [\E_0],
$$
respectively. 
\begin{theorem} \label{BCSUq2main}
Let $ q \in (-1,1) \setminus \{0\} $. Then $ \mathbb{C} $ is a retract of $ C(G_q/T) $ in $ KK^{\DD(G_q)} $. 
More precisely, we have $ \beta_q \circ \alpha_q = \id $ in $ KK^{\DD(G_q)}(\mathbb{C} \oplus \mathbb{C}, \mathbb{C} \oplus \mathbb{C}) $. 
\end{theorem} 
\proof In order to prove the assertion we have to compute the Kasparov products $ [\E_0] \circ [D] $ and $ [\E_{\pm 1}] \circ [D] $. \\
Let us first consider $ [\E_0] \circ [D] $. This class is obtained from the $ \DD(G_q) $-equivariant Fredholm module 
$ D $ by forgetting the left action of $ C(G_q/T) $. As already mentioned in the proof of proposition \ref{DiracYD}, the operator $ F $ 
intertwines the representations of $ C(G_q) $ on $ \H_1 $ and $ \H_{-1} $ induced from the $ \DD(G_q) $-Hilbert space structure.
This can be read off from the fact that the 
coefficients $ x_k(1,l,i,j) $ for $ x = a,b,c,d $ and $ k = -1,0,1 $ are symmetric in the variable $ j $. It follows 
that the resulting $ \DD(G_q) $-equivariant Kasparov $ \mathbb{C} $-$ \mathbb{C} $-module is degenerate, 
and hence $ [\E_0] \circ [D] = 0 $ in $ KK^{\DD(G_q)}(\mathbb{C}, \mathbb{C}) $. \\
Let us now study $ [\E_{-1}] \circ [D] $. The underlying graded $ \DD(G_q) $-Hilbert space of this Kasparov module is 
$$ 
\H = \H_0 \oplus \H_{-2} = L^2(\E_0) \oplus L^2(\E_{-2}), 
$$
and the corresponding action of $ C(G_q) $ is given by $ \omega $ as defined in lemma \ref{lemma0}. 
We choose the standard orthonormal basis vectors
$ e^{(l)}_{i, 0} $ for $ \H_0 $ and $ e^{(l)}_{i, -1} $ for $ \H_{-2} $. The operator 
$$
F = 
\begin{pmatrix}
0 & F_- \\
F_+ & 0 
\end{pmatrix}
$$
is determined by 
$$
F_+(e^{(l)}_{i, 0}) = 
\begin{cases}
e^{(l)}_{i, -1} & l > 0 \\
0 & l = i = 0,
\end{cases}
\qquad 
F_-(e^{(l)}_{i, -1}) = e^{(l)}_{i, 0}.
$$ 
By construction, this operator is $ G_q $-equivariant, but $ F $ does not commute with the action of 
the discrete part of $ \DD(G_q) $. \\
Let us construct a $ \DD(G_q) $-equivariant Kasparov $ \mathbb{C} $-$ C[0,1] $-module 
as follows. As underlying graded $ G_q $-equivariant Hilbert $ C[0,1] $-module we take the constant field 
$$ 
\H \otimes C[0,1] = (\H_0 \otimes C[0,1]) \oplus (\H_{-2} \otimes C[0,1])
$$ 
of $ G_q $-Hilbert spaces. It follows from lemma \ref{lemma1} that 
$$
\Omega(x)(\xi)(t) = \omega_t(x) \xi(t)
$$
defines a $ * $-homomorphism $ \Omega: C(G_q) \rightarrow \LH(\H_0 \otimes C[0,1]) $, 
and the corresponding coaction of $ C^*(G_q) $ turns the even part $ \H_0 \otimes C[0,1] $ into 
a $ \DD(G_q) $-equivariant Hilbert $ C[0,1] $-module. 
In odd degree we consider the constant $ \DD(G_q) $-Hilbert module structure induced from $ \H_{-2} $. 
The left action of $ \mathbb{C} $ on $ \H \otimes C[0,1] $ is given by multiples of the 
identity operator, and as a final ingredient we take the constant operator $ F \otimes 1 $ 
on $ \H \otimes C[0,1] $. By construction, this operator is $ G_q $-equivariant, and it follows from lemma \ref{lemma2} 
that $ F \otimes 1 $ commutes with the action of $ C(G_q) $ up to compact operators. Hence 
we have indeed defined a Kasparov module. \\
Evaluation of this Kasparov module at $ t = 1 $ yields $ [\E_{-1}] \circ [D] $. 
The evaluation at $ t = 0 $ agrees with the cycle defining $ [\E_{-1}] \circ [D] $ except that 
the action of $ C(G_q) $ on $ \H_0 $ is given by $ \omega_0 $ instead of $ \omega_1 $. 
We decompose 
$$ 
\H_0 =  \mathbb{C} \oplus \mathbb{C}^\bot 
$$ 
as a direct sum of the one-dimensional subspace $ \mathbb{C} $ spanned by $ e^{(0)}_{0,0} $ and 
its orthogonal complement $ \mathbb{C}^\bot $. Inspection of the explicit formulas shows that $ \omega_0 $ preserves this 
decomposition and implements the trivial representation $ \epsilon: C(G_q) \rightarrow \mathbb{C} = \LH(\mathbb{C}) $ 
on the first component. It follows that the Kasparov module decomposes as a direct sum of the trivial module 
$ \mathbb{C} $ representing the identity and its orthogonal complement which we denote by $ R $. 
That is, the underlying graded $ \DD(G_q) $-Hilbert space of $ R $ is $ \mathbb{C}^\bot \oplus \H_{-2} $, the 
representation of $ \mathbb{C} $ is given by multiples of the identity, and the operator $ F $ defines a 
$ G_q $-equivariant isomorphism between $ \mathbb{C}^\bot $ and $ \H_{-2} $. \\
For $ q > 0 $ we see from lemma \ref{lemma3} that $ F $ intertwines the representation $ \omega_0 $ on 
$ \mathbb{C}^\bot $ with the representation $ \omega $ on $ \H_{-2} $. It follows that $ R $ is degenerate, 
and hence $ [\E_{-1}] \circ [D] = \id $ in $ KK^{\DD(G_q)}(\mathbb{C}, \mathbb{C}) $ 
in this case. \\
For $ q < 0 $ the module $ R $ fails to be degenerate since the representations of $ C(G_q) $ 
do not match. Let $ F $ be the operator on $ \H_{-2} \oplus \H_{-2} $ given by
$$ 
F = 
\begin{pmatrix}
0 & 1 \\
1 & 0 
\end{pmatrix}.
$$
We define a homotopy $ T $ by considering the constant field of graded $ G_q $-Hilbert spaces 
$$
\H \otimes C[0,1] = (\H_{-2} \oplus \H_{-2}) \otimes C[0,1] \oplus (\H_{-2} \oplus \H_{-2}) \otimes C[0,1],
$$
the representation $ \phi: \mathbb{C} \rightarrow \LH(\H \otimes C[0,1]) $ given by
$$ 
\phi = 
\begin{pmatrix}
\phi_+ & 0 \\
0 & \phi_-
\end{pmatrix}
$$
where
$$ 
\phi_\pm(\lambda)
\begin{pmatrix}
\xi \\
\eta
\end{pmatrix}
= 
\begin{pmatrix}
\lambda \, \xi \\
0
\end{pmatrix}
$$
and the constant operator $ (F \oplus F) \otimes 1 $. \\
Let us use the canonical identification of $ \mathbb{C}^\bot $ with $ \H_{-2} $ 
in order to view $ \omega_0 $ as a representation on $ \H_{-2} $. 
We define a $ \DD(G_q) $-Hilbert module structure on the even part $ (\H_{-2} \oplus \H_{-2}) \otimes C[0,1] $ of $ T $ 
by the action 
\begin{align*}
\Theta_+(x)
\begin{pmatrix}
\xi \\
\eta
\end{pmatrix}
(t)
&= 
U(t) 
\begin{pmatrix}
\omega_0(x) & 0 \\
0 & \omega(x)
\end{pmatrix}
U^{-1}(t)
\begin{pmatrix}
\xi(t) \\
\eta(t)
\end{pmatrix} 
\end{align*}
of $ C(G_q) $ where $ U(t) $ is the rotation matrix
$$
U(t) = \begin{pmatrix}
\cos(\pi t/2) & \sin(\pi t/2) \\
- \sin(\pi t/2) & \cos(\pi t/2)
\end{pmatrix}.
$$
In odd degree we consider the constant $ \DD(G_q) $-Hilbert module structure 
$$ 
\Theta_-(x)
\begin{pmatrix}
\xi \\
\eta
\end{pmatrix} 
= 
\begin{pmatrix}
\omega(x) & 0 \\
0 & \omega(x) 
\end{pmatrix}
\begin{pmatrix}
\xi \\
\eta
\end{pmatrix} 
= 
\begin{pmatrix}
\omega(x)(\xi) \\
\omega(x)(\eta)
\end{pmatrix}. 
$$
Writing $ \Theta = \Theta_+ \oplus \Theta_- $ we see from lemma \ref{lemma2} and lemma \ref{lemma3} that the commutators 
$ [(F \oplus F) \otimes 1, \Theta(x)] $ 
are compact for all $ x \in C(G_q) $. It follows that $ T $ defines a $ \DD(G_q) $-equivariant 
Kasparov $ \mathbb{C} $-$ C[0,1] $-module. \\
The evaluation of $ T $ at $ t = 0 $ identifies with the sum of $ R $ and a degenerate module, and evaluation at $ t = 1 $ 
yields a degenerate module. We conclude that the Kasparov module $ R $ is homotopic to zero, and hence 
$ [\E_{-1}] \circ [D] = \id $ in $ KK^{\DD(G_q)}(\mathbb{C}, \mathbb{C}) $. \\
In a similar way one proves the relation $ (-[\E_1]) \circ [D] = \id $ in $ KK^{\DD(G_q)}(\mathbb{C}, \mathbb{C}) $ for all 
$ q \in (-1,1) \setminus \{0\} $. 
The calculations are analogous and will be omitted. \qed \\
As a corollary of theorem \ref{BCSUq2main} we obtain the following result. 
\begin{theorem} \label{PDPodles}
Let $ q \in (0,1) $. The standard Podle\'s sphere $ C(G_q/T) $ is isomorphic to 
$ \mathbb{C} \oplus \mathbb{C} $ in $ KK^{\DD(G_q)} $. 
\end{theorem} 
\proof According to theorem 6.7 in \cite{NVpoincare} the elements $ \alpha_q $ and $ \beta_q $ considered above 
satisfy $ \alpha_q \circ \beta_q = \id $ in $ KK^{\DD(G_q)}(C(G_q/T),C(G_q/T)) $. 
Hence due to theorem \ref{BCSUq2main} these elements induce inverse isomorphisms in $ KK^{\DD(G_q)} $ as desired. \qed \\
We have not checked wether the assertion of theorem \ref{PDPodles} holds for $ q < 0 $ as well. For our purposes 
the following result is sufficient, compare \cite{NVpoincare}.
\begin{prop} \label{alphabetaid}
Let $ q \in (-1,1) \setminus \{0\} $. The standard Podle\'s sphere $ C(G_q/T) $ is isomorphic to $ \mathbb{C} \oplus \mathbb{C} $ in $ KK^{G_q} $.
\end{prop}
\proof For $ q > 0 $ this is an immediate consequence of theorem \ref{PDPodles}. Hence it remains 
to treat the case of negative $ q $, and clearly it suffices to show $ \alpha_q \circ \beta_q = \id $ in $ KK^{G_q}(C(G_q/T), C(G_q/T)) $ for 
all $ q \in (-1,1) \setminus \{0\} $. \\
Note first that the definition of the $ C^* $-algebras $ C(G_q) $ and $ C(G_q/T) $ 
makes sense also for $ q = 0 $. Moreover, the algebras $ C(G_q/T) $ assemble into a $ T $-equivariant continuous 
field over $ [q,1] $ for all $ q \in (-1,1] $, see 
\cite{Blanchard}, \cite{Nagydeformation}. Let us write $ C({\bf G}/T) $ for the corresponding algebra of sections. 
The elements $ \alpha_t $ yield a class $ \alpha \in KK^T(C({\bf G}/T), C[q,1]) $, 
in particular, we have a well-defined element $ \alpha_0 $ for $ t = 0 $. 
Similarly, the elements $ \beta_t $ determine a class $ \beta \in KK^T(C[q,1], C({\bf G}/T)) $. 
Using a comparison argument as in \cite{NVpoincare}, the claim follows from the induction isomorphism 
$$
KK^{G_q}(C(G_q/T), C(G_q/T)) \cong KK^T(C(G_q/T), \mathbb{C}) 
$$
and the fact that $ \alpha_1 \circ \beta_1 = \id \in KK^{G_1}(C(G_1/T), C(G_1/T)) $. \qed

\section{The Baum-Connes conjecture for torsion-free discrete quantum groups} \label{secbc}

In this section we recall the formulation of the Baum-Connes conjecture for torsion-free discrete quantum groups
proposed by Meyer \cite{Meyerhomalg2}. This involves some general concepts from homological algebra 
in triangulated categories. For more detailed information we refer to \cite{Neeman}, \cite{MNhomalg1}, \cite{Meyerhomalg2}, 
\cite{MNtriangulated}. \\
Let $ G $ be a torsion-free discrete quantum group in the sense of definition \ref{deftorsion}. 
The equivariant Kasparov category $ KK^G $ has as objects all separable
$ G $-$ C^* $-algebras, and $ KK^G(A,B) $ as the set of morphisms between two objects $ A $ and $ B $.
Composition of morphisms is given by the Kasparov product. The category $ KK^G $ is triangulated with 
translation automorphism $ \Sigma: KK^G \rightarrow KK^G $ given by 
the suspension $ \Sigma A = C_0(\mathbb{R}, A) $ of a $ G $-$ C^* $-algebra $ A $.
Every $ G $-equivariant $ * $-homomorphism $ f: A \rightarrow B $ induces a diagram of the form
$$
\xymatrix{
\Sigma B  \;\; \ar@{->}[r] & C_f \ar@{->}[r] & A \ar@{->}[r]^f & B
}
$$
where $ C_f $ denotes the mapping cone of $ f $. Such diagrams are called mapping cone triangles. 
By definition, an exact triangle is a diagram in $ KK^G $ of the form $ \Sigma Q \rightarrow K \rightarrow E \rightarrow Q $
which is isomorphic to a mapping cone triangle. \\
Associated with the inclusion of the trivial quantum subgroup $ E $ in $ G $ we have the obvious restriction functor 
$ \res^G_E: KK^G \rightarrow KK^E = KK $ 
and an induction functor $ \ind_E^G: KK \rightarrow KK^G $. Explicitly, 
$ \ind_E^G(A) = C_0(G) \otimes A $ for $ A \in KK $ with action of $ G $ given by translation on the copy of $ C_0(G) $. \\
We consider the following full subcategories of $ KK^G $, 
\begin{align*}
\CC_G &= \{A \in KK^G|\res^G_E(A) = 0 \in KK \} \\
\CI_G &= \{\ind_E^G(A)| A \in KK \}, 
\end{align*}
and refer to their objects as compactly contractible and compactly induced 
$ G $-$ C^* $-algebras, respectively. Since $ G $ is torsion-free, it suffices to consider the trivial quantum 
subgroup in the definition of these categories. 
If there is no risk of confusion we will write 
$ \CC $ and $ \CI $ instead of $ \CC_G $ and $ \CI_G $. \\
The subcategory $ \CC $ is localising, and we denote by $ \bra \CI \ket $ the localising 
subcategory generated by $ \CI $. It follows from theorem 3.21 in \cite{Meyerhomalg2} that the pair of localising subcategories 
$ (\bra \CI\ket, \CC) $ in $ KK^G $ is complementary. That is, $ KK^G(P,N) = 0 $ for 
all $ P \in \bra \CI \ket $ and $ N \in \CC $, and every object $ A \in KK^G $ fits into an exact triangle 
$$
\xymatrix{
\Sigma N \; \ar@{->}[r] & \tilde{A} \ar@{->}[r] & A \ar@{->}[r] & N
}
$$
with $ \tilde{A} \in \bra \CI \ket $ and $ N \in \CC $. 
Such a triangle is uniquely determined up to isomorphism and depends functorially on $ A $. 
We will call the morphism $ \tilde{A} \rightarrow A $ a Dirac element for $ A $. \\
The localisation $ \mathbb{L}F $ of a homological functor $ F $ on $ KK^G $ 
at $ \CC $ is given by 
$$
\mathbb{L}F(A) = F(\tilde{A}) 
$$
where $ \tilde{A} \rightarrow A $ is a Dirac element for $ A $. 
By construction, there is an obvious map $ \mathbb{L}F(A) \rightarrow F(A) $ for 
every $ A \in KK^G $. \\
In the sequel we write 
$ G \ltimes_\max A $ and $ G \ltimes_\red A $ for the full and reduced crossed products of $ A $ by $ G $. 
Let us remark that in \cite{NVpoincare} these algebras are denoted by $ C^*_\max(G)^\cop \ltimes_\max A $ and $ C^*_\red(G)^\cop \ltimes_\red A $, 
respectively. 
\begin{definition}
Let $ G $ be a torsion-free discrete quantum group and consider the functor $ F(A) = K_*(G \ltimes_\red A) $ on $ KK^G $. The 
Baum-Connes assembly map for $ G $ with coefficients in $ A $ is the map 
$$ 
\mu_A: \mathbb{L}F(A) \rightarrow F(A). 
$$
We say that $ G $ satisfies the Baum-Connes conjecture with coefficients in $ A $ if $ \mu_A $ is an isomorphism. 
We say that $ G $ satisfies the strong Baum-Connes conjecture if $ \bra \CI \ket = KK^G $. 
\end{definition}
Observe that the strong Baum-Connes conjecture implies the Baum-Connes conjecture with arbitrary coefficients. 
Indeed, for $ A \in \bra \CI \ket $ the assembly map $ \mu_A $ is clearly an isomorphism. \\
By the work of Meyer and Nest \cite{MNtriangulated}, the above terminology is consistent with the classical definitions 
in the case that $ G $ is a torsion-free discrete group. 
The strong Baum-Connes conjecture amounts to the assertion that $ G $ has a $ \gamma $-element 
and $ \gamma = 1 $ in this case. \\
In section \ref{secapp} we will need further considerations from \cite{MNhomalg1}, \cite{Meyerhomalg2} 
relying on the notion of a homological ideal in a triangulated category. Let us briefly discuss the relevant material 
adapted to our specific situation. \\
We denote by $ \mathfrak{J} $ the homological ideal in $ KK^G $ consisting of all $ f \in KK^G(A,B) $ such that 
$ \res^G_E(f) = 0 \in KK(A,B) $. By definition, $ \mathfrak{J} $ is the kernel of the exact functor 
$ \res^G_E: KK^G \rightarrow KK $. The ideal $ \mathfrak{J} $ is compatible with countable direct sums and 
has enough projective objects. The $ \mathfrak{J} $-projective objects in $ KK^G $ are precisely the retracts 
of compactly induced $ G $-$ C^* $-algebras. \\
A chain complex 
\begin{equation*}
\xymatrix{
\cdots \ar@{->}[r] & C_{n + 1} \ar@{->}[r]^{d_{n + 1}} & C_n \ar@{->}[r]^{d_n} & C_{n - 1} \ar@{->}[r] & \cdots 
}
\end{equation*}
in $ KK^G $ is $ \mathfrak{J} $-exact if 
\begin{equation*}
\xymatrix{
\cdots \ar@{->}[r] & KK(A, C_{n + 1}) \ar@{->}[r]^{\;\;\;(d_{n + 1})_*} & KK(A, C_n) \ar@{->}[r]^{(d_n)_*} & KK(A, C_{n - 1}) \ar@{->}[r] & \cdots 
}
\end{equation*}
is exact for every $ A \in KK $. \\
A $ \mathfrak{J} $-projective resolution of $ A \in KK^G $ is a chain complex 
\begin{equation*}
\xymatrix{
\cdots \ar@{->}[r] & P_{n + 1} \ar@{->}[r]^{d_{n + 1}} & P_n \ar@{->}[r]^{d_n} & P_{n - 1} \ar@{->}[r] & \cdots \ar@{->}[r]^{d_2} & P_1 \ar@{->}[r]^{d_1} & P_0
}
\end{equation*}
of $ \mathfrak{J} $-projective objects in $ KK^G $, augmented by a map $ P_0 \rightarrow A $ such that 
the augmented chain complex is $ \mathfrak{J} $-exact. \\
For our purposes it is important that a $ \mathfrak{J} $-projective resolution of $ A \in KK^G $ can be used to 
construct a Dirac element $ \tilde{A} \rightarrow A $. In general, this construction leads to a spectral sequence 
computing the derived functor $ \mathbb{L}F(A) $. In the specific case of free orthogonal quantum groups 
that we are interested in, the spectral sequence reduces to a short exact sequence. This short exact sequence 
will be discussed in section \ref{secapp} in connection with our $ K $-theory computations.

\section{The Baum-Connes conjecture for the dual of $ SU_q(2) $} \label{secbcsuq2}

In this section we show that the dual of $ SU_q(2) $ satisfies the strong Baum-Connes conjecture. 
We work within the general setup explained in the previous section, taking into account proposition \ref{suq2torsion}
which asserts that the dual of $ SU_q(2) $ is torsion-free. Let us remark that the strong Baum-Connes conjecture for the 
dual of the classical group $ SU(2) $ is a special case of the results in \cite{MNcompact}.
\begin{theorem}\label{BCsuq2}
Let $ q \in (-1,1) \setminus \{0\} $. The dual discrete quantum group of $ SU_q(2) $ satisfies the strong Baum-Connes conjecture. 
\end{theorem}
\proof In the sequel we write $ G = SU_q(2) $. Due to Baaj-Skandalis 
duality it suffices to prove that every $ G $-$ C^* $-algebra is contained in the localising 
subcategory $ \T $ of $ KK^G $ generated by all trivial $ G $-$ C^* $-algebras. \\
Let $ A $ be a $ G $-$ C^* $-algebra. Theorem \ref{BCSUq2main} implies that $ A $ is a retract of $ C(G/T) \twisted_G A $ in $ KK^G $, 
and according to theorem 3.6 in \cite{NVpoincare} there is a natural $ G $-equivariant 
isomorphism $ C(G/T) \twisted_G A \cong \ind_T^G \res_T^G(A) $. \\
Since $ \hat{T} = \mathbb{Z} $ is a torsion-free discrete abelian group the strong
Baum-Connes conjecture holds for $ \hat{T} $. More precisely, the trivial $ \hat{T} $-$ C^* $-algebra $ \mathbb{C} $ is contained in 
the localising subcategory $ \bra C_0(\hat{T}) \ket $ of $ KK^{\hat{T}} $ generated by $ C_0(\hat{T}) $. 
Next observe that there is a $ \hat{T} $-equivariant $ * $-isomorphism 
$ C_0(\hat{T}) \otimes B \cong C_0(\hat{T}) \otimes B_\tau $ for every $ \hat{T} $-$ C^* $-algebra $ B $ where $ B_\tau $ denotes $ B $ 
with the trivial $ \hat{T} $-action. It follows that 
\begin{align*} 
T \ltimes \res^G_T(A) \cong \mathbb{C} \otimes T \ltimes \res^G_T(A) \in \bra C_0(\hat{T}) &\otimes T \ltimes \res^G_T(A) \ket \\
&= \bra C_0(\hat{T}) \otimes (T \ltimes \res^G_T(A))_\tau \ket
\end{align*}
in $ KK^{\hat{T}} $. According to Baaj-Skandalis duality, this implies
\begin{align*}
\res^G_T(A) \cong \hat{T} \ltimes T \ltimes \res^G_T(A) \in
\bra \hat{T} \ltimes &(C_0(\hat{T}) \otimes (T \ltimes \res^G_T(A))_\tau) \ket \\
&= \bra (T \ltimes \res^G_T(A))_\tau \ket  
\end{align*}
in $ KK^T $. Using proposition \ref{alphabetaid} we thus obtain 
$$
\ind_T^G \res^G_T(A) \in \bra C(G/T) \otimes (T \ltimes \res^G_T(A))_\tau \ket \subset \bra (T \ltimes \res^G_T(A))_\tau \ket \subset \T 
$$
in $ KK^G $ since the induction functor $ \ind_T^G $ is triangulated. \\
Combining the above considerations shows $ A \in \T $ and finishes the proof. \qed \\
Starting from theorem \ref{BCsuq2} it is easy to calculate the $ K $-groups of 
$ C(SU_q(2)) $ and $ C(SO_q(3)) $. We shall not present these computations here. 
 
\section{Free orthogonal quantum groups and monoidal equivalence} \label{secfreeqg}

In this section we review the definition of free orthogonal quantum groups and discuss 
the concept of monoidal equivalence for compact quantum groups. \\
We begin with the definition of free orthogonal quantum groups. These quantum groups were introduced by Wang and van Daele \cite{Wang}, \cite{vDW}. 
As usual, for a matrix $ u = (u_{ij}) $ of elements in a $ * $-algebra we shall write
$ \overline{u} $ and $ u^t $ for its conjugate and transposed matrices, respectively. 
That is, we have $ (\overline{u})_{ij} = u_{ij}^* $ and $ (u^t)_{ij} = u_{ji} $ for the corresponding matrix entries. 
\begin{definition} \label{deffo}
Let $ Q \in GL_n(\mathbb{C}) $ such that $ Q \overline{Q} = \pm 1 $. The group $ C^* $-algebra $ C^*_\max(\mathbb{F}O(Q)) $ 
of the free orthogonal quantum group $ \mathbb{F}O(Q) $ is the universal $ C^* $-algebra with generators 
$ u_{ij}, 1 \leq i,j \leq n $ such that the resulting matrix $ u $ is unitary and the relation
$ u = Q \overline{u} Q^{-1} $ holds. 
\end{definition} 
In definition \ref{deffo} we basically adopt the conventions in \cite{Banicaunitary}. However, we write $ Q $ instead 
of $ F $ for the parameter matrix, and we deviate from the standard notation $ A_o(Q) = C^*_\max(\mathbb{F}O(Q)) $. The latter
is motivated from the fact that we shall view this $ C^* $-algebra as the group $ C^* $-algebra of a discrete quantum group. \\
It is well-known that $ C(SU_q(2)) $ can be written as the group $ C^* $-algebra of a free orthogonal quantum group for 
appropriate $ Q \in GL_2(\mathbb{C}) $. Moreover, the free quantum groups $ \mathbb{F}O(Q) $ for $ Q \in GL_2(\mathbb{C}) $ 
exhaust up to isomorphism precisely the duals of $ SU_q(2) $ for $ q \in [-1,1] \setminus \{0\} $. \\
The quantum groups $ \mathbb{F}O(Q) $ for higher dimensional matrices $ Q $ are still closely related to 
quantum $ SU(2) $. In order to explain this, we shall discuss the notion of monoidal equivalence for compact quantum groups 
introduced by Bichon, de Rijdt and Vaes \cite{BdRV}. For the algebraic aspects of monoidal equivalences and Hopf-Galois theory 
we refer to \cite{Schauenburg}. \\
As in section \ref{secspec} we write $ \Rep(G) $ for the $ C^* $-tensor category 
of finite dimensional representations of a compact quantum group $ G $. Recall that the objects in $ \Rep(G) $ 
are the finite dimensional representations of $ G $, and the morphism sets consist
of all intertwining operators. 
\begin{definition} \label{defme}
Two compact quantum groups $ G $ and $ H $ are called monoidally equivalent if the 
representation categories $ \Rep(G) $ and $ \Rep(H) $ are equivalent as $ C^* $-tensor categories. 
\end{definition}
We point out that in definition \ref{defme} the representation categories are only required to be equivalent as abstract $ C^* $-tensor 
categories. In fact, by the Tannaka-Krein reconstruction theorem \cite{Woronowiczsun}, a compact quantum group $ G $ is determined 
up to isomorphism by the $ C^* $-tensor category $ \Rep(G) $ together with its canonical 
fiber functor into the category of Hilbert spaces. \\
Let $ H $ be a compact quantum group. An algebraic coaction $ \lambda: \P \rightarrow \mathbb{C}[H] \odot \P $ on the unital $ * $-algebra $ \P $ is 
called ergodic if the invariant subalgebra $ \mathbb{C} \Box_{\mathbb{C}[H]} \P \subset \P $ is equal to $ \mathbb{C} $. 
We say that $ \P $ is a left Galois object if $ \lambda $ is ergodic and the Galois map 
$ \gamma_\P: \P \odot \P \rightarrow \mathbb{C}[H] \odot \P $ given by 
$$ 
\gamma_\P(x \odot y) = \lambda(x) (1 \odot y) 
$$ 
is a linear isomorphism. 
Similarly one defines ergodicity for right coactions and the notion of a right Galois object. \\
Like in Morita theory, it is important that monoidal equivalences can be implemented concretely.  
\begin{definition}
Let $ G $ and $ H $ be compact quantum groups. A bi-Galois object for $ G $ and $ H $ is a unital $ * $-algebra $ \P $ which is both a 
left $ \mathbb{C}[H] $-Galois object and a right $ \mathbb{C}[G] $-Galois object, such that the 
corresponding coactions turn $ \P $ into a $ \mathbb{C}[H] $-$ \mathbb{C}[G] $-bicomodule. 
\end{definition} 
A linear functional $ \omega $ on a unital $ * $-algebra $ \P $ is called a state if $ \omega(x^*x) \geq 0 $ 
for all $ x \in \P $ and $ \omega(1) = 1 $. A state $ \omega $ is said to be faithful if $ \omega(x^*x) = 0 $ implies $ x = 0 $. 
If $ \P $ is in addition equipped with a coaction $ \lambda: \P \rightarrow \mathbb{C}[H] \odot \P $ then $ \omega $ is called 
invariant if $ (\id \odot \omega)\lambda(x) = \omega(x) 1 $ for all $ x \in \P $. \\
The following result is proved in \cite{BdRV}. 
\begin{theorem} \label{BdRVmain} 
Let $ G $ and $ H $ be monoidally equivalent compact quantum groups. 
Then there exists a bi-Galois object $ \P $ for $ G $ and $ H $ such that 
$$ 
\F(\H) = \P \Box_{\mathbb{C}[G]} \H 
$$ 
defines a monoidal equivalence $ \F: \Rep(G) \rightarrow \Rep(H) $, and there exists a canonical faithful state $ \omega $ on 
$ \P $ which is left and right invariant with respect to the coactions of $ \mathbb{C}[G] $ and 
$ \mathbb{C}[H] $, respectively.  
\end{theorem}
As a first application of the concept of monoidal equivalence let us record the following fact. 
\begin{prop} \label{MEtorsion} 
Let $ G $ and $ H $ be discrete quantum groups with monoidally equivalent duals. Then $ G $ is torsion-free 
iff $ H $ is torsion-free. 
\end{prop}
\proof As explained in \cite{dRV}, actions of monoidally equivalent compact quantum groups 
on finite dimensional $ C^* $-algebras are in a bijective correspondence. We will discuss this more generally,  
for arbitrary $ C^* $-algebras and on the level of equivariant $ KK $-theory, in section \ref{secmon}. Under this 
correspondence, actions associated to representations of $ G $ correspond to actions associated to representations of $ H $. 
This immediately yields the claim. \qed \\
In the sequel we will make use of the following crucial result from \cite{BdRV}, which in turn relies on the 
fundamental work of Banica \cite{Banicafo}, \cite{Banicaunitary}. 
\begin{theorem} \label{MEfo}
Let $ Q_j \in GL_{n_j}(\mathbb{C}) $ such that $ Q_j \overline{Q_j} = \pm 1 $ for $ j = 1,2 $. Then the dual of 
$ \mathbb{F}O(Q_1) $ is monoidally equivalent to the dual of $ \mathbb{F}O(Q_2) $ iff 
$ Q_1 \overline{Q_1} $ and $ Q_2 \overline{Q_2} $ have the same sign and 
$$
\tr(Q_1^* Q_1) = \tr(Q_2^* Q_2). 
$$
In particular, for any $ Q \in GL_n(\mathbb{C}) $ such that $ Q \overline{Q} = \pm 1 $, 
the dual of $ \mathbb{F}O(Q) $ is monoidally equivalent to $ SU_q(2) $ 
for a unique $ q \in [-1,1] \setminus \{0\} $.  
\end{theorem}
Theorem \ref{MEfo} implies in particular that the dual of $ \mathbb{F}O(Q) $ for $ Q \in GL_n(\mathbb{C}) $ and $ n > 2 $
is not monoidally equivalent to $ SU_{\pm 1}(2) $. With this in mind 
we obtain the following consequence of proposition \ref{suq2torsion} and proposition \ref{MEtorsion}. 
\begin{cor} 
Let $ Q \in GL_n(\mathbb{C}) $ for $ n > 2 $ such that $ Q \overline{Q} = \pm 1 $.  
Then the free orthogonal quantum group $ \mathbb{F}O(Q) $ is torsion-free. 
\end{cor}

\section{Monoidal equivalence and equivariant $ KK $-theory} \label{secmon}

Extending considerations in \cite{dRV}, we discuss in this section the correspondence for actions of monoidally equivalent 
compact quantum groups. We show in particular that the strong Baum-Connes property
for torsion-free quantum groups is invariant under monoidal equivalence. \\
Let $ G $ and $ H $ be monoidally equivalent compact quantum groups 
and let $ \P $ be the bi-Galois object for $ G $ and $ H $ as in theorem \ref{BdRVmain}. 
Moreover let $ A $ be a $ G $-$ C^* $-algebra, and recall from section \ref{secspec} that we write $ \S(A) $ for the dense 
spectral $ * $-subalgebra of $ A $. The algebraic cotensor product 
$ \F(A) = \P \Box_{\mathbb{C}[G]} \S(A) \subset \P \odot \S(A) $ is again a $ * $-algebra 
and carries an algebraic coaction $ \lambda: \F(\A) \rightarrow \mathbb{C}[H] \odot \F(A) $ inherited from $ \P $. \\
Consider the $ C^* $-algebra $ P \otimes A $ where $ P $ denotes the minimal completion of $ \P $, that 
is, the $ C^* $-algebra generated by $ \P $ in the GNS-representation of the invariant state $ \omega $. 
The left coaction on $ P $ turns $ P \otimes A $ into an $ H $-$ C^* $-algebra.  
Let $ F(A) = P \Box_G A $ be the closure of $ \F(A) = \P \Box_{\mathbb{C}[G]} \S(A) $ inside $ P \otimes A $. By construction, 
the coaction of $ P \otimes A $ maps $ \F(A) $ into $ C^\red(H) \otimes F(A) $. 
In this way we obtain a coaction on $ F(A) $ which turns $ F(A) $ into an $ H $-$ C^* $-algebra. \\
If $ f: A \rightarrow B $ is a $ G $-equivariant $ * $-homomorphism then $ \id \otimes f: P \otimes A \rightarrow P \otimes B $ 
induces an $ H $-equivariant $ * $-homomorphism $ \id \Box_G f: P \Box_G A \rightarrow P \Box_G B $. 
Consequently, we obtain a functor $ F: G \Alg \rightarrow H \Alg $ by setting $ F(A) = P \Box_G A $ on objects and $ F(f) = \id \Box_G f $
on morphisms. Here $ G\Alg $ and $ H \Alg $ denote the categories of separable $ G $-$ C^*$-algebras 
and $ H $-$ C^* $-algebras, respectively. Note that a trivial $ G $-$ C^* $-algebra $ A $ is mapped to the trivial $ H $-$ C^* $-algebra 
$ F(A) \cong A $ under the functor $ F $. Moreover $ F(A \oplus B) \cong F(A) \oplus F(B) $ for all 
$ G $-$ C^* $-algebras $ A $ and $ B $.\\
By symmetry, we have the dual Galois object $ \Q $ for $ H $ and $ G $ and a corresponding functor 
$ H \Alg \rightarrow G \Alg $. This functor sends an $ H $-$ C^* $-algebra $ B $ to the $ G $-$ C^* $-algebra $ Q \Box_H B $. 
Here $ Q $ denotes the $ C^* $-algebra generated by $ \Q $ in the GNS-representation associated to its natural invariant state $ \eta $. 
\begin{prop} \label{equivlemma1}
For every $ G $-$ C^* $-algebra $ A $ there is a natural isomorphism 
$$
Q \Box_H P \Box_G A \cong A 
$$
of $ G $-$ C^* $-algebras. 
\end{prop}
\proof Consider first the case $ A = C^\red(G) $. In this case we have a canonical isomorphism $ P \Box_G A \cong P $. 
By construction, $ Q \Box_H P \subset Q \otimes P $ is the closure of $ \Q \Box_{\mathbb{C}[H]} \P \cong \mathbb{C}[G] $. 
Since $ \eta $ and $ \omega $ are faithful states on $ Q $ and $ P $, respectively, the state $ \eta \otimes \omega $ 
is faithful on $ Q \otimes P $, and hence also on $ Q \Box_H P $. 
Moreover $ \eta \otimes \omega $ is both left and right invariant with respect to the natural coactions of $ C^\red(G) $. 
From this we conclude that the above inclusion $ \mathbb{C}[G] \rightarrow Q \otimes P $ induces an equivariant 
$ * $-isomorphism $ C^\red(G) \cong Q \Box_H P $. \\
Now let $ A $ be an arbitrary $ G $-$ C^* $-algebra. Using the previous discussion we 
obtain a well-defined injective $ * $-homomorphism $ \alpha: A \rightarrow Q \otimes P \otimes A $ 
by applying the coaction followed with the isomorphism $ C^\red(G) \cong Q \Box_H P \subset Q \otimes P $. 
Due to associativity of the cotensor product the coaction $ \S(A) \rightarrow \mathbb{C}[G] \odot \S(A) $ induces 
an isomorphism
$$
\S(A) \cong \mathbb{C}[G] \Box_{\mathbb{C}[G]} \S(A) \cong \Q \Box_{\mathbb{C}[H]} \P \Box_{\mathbb{C}[G]} \S(A). 
$$
Now let $ \pi \in \hat{H} $. The spectral subspace $ _\pi P $ is a finite dimensional right $ \mathbb{C}[G] $-comodule, 
and we observe that $ _\pi(P \otimes A) = (_\pi P) \odot A $. In fact, 
we have $ _\pi(P \otimes A) = (p_\pi \otimes \id)(P \otimes A) \subset (_\pi P) \odot A $ where $ p_\pi: P \rightarrow\, _\pi P $ 
is the projection operator, and the reverse inclusion is obvious. This implies
$ _\pi(P \Box_G A) = (_\pi P) \Box_{\mathbb{C}[G]} \S(A) $ 
and hence $ \S(P \Box_G A) = \P \Box_{\mathbb{C}[G]} \S(A) $ for the spectral subalgebras. Using a symmetric 
argument for $ Q $ we conclude
$$
\S(Q \Box_H P \Box_G A) = \Q \Box_{\mathbb{C}[H]} \S(P \Box_G A) = \Q \Box_{\mathbb{C}[H]} \P \Box_{\mathbb{C}[G]} \S(A). 
$$
It follows that the image of $ \S(A) $ in $ Q \Box_H P \Box_G A $ under the map $ \alpha $ is dense. 
Hence $ \alpha $ induces an equivariant $ * $-isomorphism $  A \cong Q \Box_H P \Box_G A $ as desired. \qed \\
We have thus proved the following theorem. 
\begin{theorem} \label{MEalg}
Let $ H $ and $ G $ be monoidally equivalent compact quantum groups. Then the categories $ G\Alg $ and $ H\Alg $ 
are equivalent. 
\end{theorem}
Our next aim is to extend the equivalence of theorem \ref{MEalg} to the level of equivariant Kasparov theory. \\
For this we need an appropriate version of the cotensor product for Hilbert modules. Let $ \E $ 
be a $ G $-Hilbert $ A $-module, and recall from section \ref{secspec} 
that $ \S(\E) $ denotes the spectral submodule of $ \E $. 
We define $ P \Box_G \E $ as the closure of $ \P \Box_{\mathbb{C}[G]} \S(\E) $ inside the Hilbert $ P \otimes A $-module $ P \otimes \E $. 
It is straightforward to check that $ P \Box_G \E $ is an $ H $-Hilbert $ P \Box_G A $-module. Moreover, the following assertion is proved in 
the same way as proposition \ref{equivlemma1}. 
\begin{prop} \label{MEisohilbert}
For every $ G $-Hilbert $ A $-module $ \E $ there is a natural isomorphism 
$$
Q \Box_H P \Box_G \E \cong \E 
$$
of $ G $-Hilbert $ A $-modules. 
\end{prop} 
If $ \E $ is a $ G $-Hilbert $ A $-module then $ \KH(\E) $ is a $ G $-$ C^* $-algebra in a natural way. 
The cotensor product constructions for $ C^* $-algebras and Hilbert modules are compatible in the following 
sense. 
\begin{prop} 
Let $ \E $ be a $ G $-Hilbert $ A $-module. Then 
$$
\KH(P \Box_G \E) \cong P \Box_G \KH(\E) 
$$
as $ H $-$ C^* $-algebras. 
\end{prop}
\proof Note that there are canonical inclusions $ P \Box_G \KH(\E) \subset P \otimes \KH(\E) \cong \KH(P \otimes \E) $ and 
$ \KH(P \Box_G \E) \subset \KH(P \otimes \E) $. These inclusions determine a homomorphism 
$ \iota_\E: \KH(P \Box_G \E) \rightarrow P \Box_G \KH(\E) $ of $ H $-$ C^* $-algebras. Using proposition \ref{MEisohilbert} 
and proposition \ref{equivlemma1} we obtain a map 
$$
P \Box_G \KH(\E) \cong P \Box_G \KH(Q \Box_H P \Box_G \E) \rightarrow P \Box_G Q \Box_H \KH(P \Box_G \E) \cong \KH(P \Box_G \E)
$$
where the middle arrow is given by $ \id \Box \iota_{P \Box_G \E} $. It is readily checked that this map is inverse to 
$ \iota_\E $.  \qed \\
Let $ \E $ and $ \F $ be $ G $-Hilbert-$ A $-modules and let $ T \in \LH(\E,\F) $ be $ G $-equivariant. Then 
$ \id \otimes T: P \otimes \E \rightarrow P \otimes \F $ 
induces an adjointable operator $ \id \Box_G T: P \Box_G \E \rightarrow P \Box_G \F $. 
If $ \phi: A \rightarrow \LH(\E) $ 
is a $ G $-equivariant $ * $-homomorphism then $ \id \otimes \phi: P \otimes A \rightarrow \LH(P \otimes \E) $
induces an $ H $-equivariant $ * $-homomorphism $ \id \Box_G \phi: P \Box_G A \rightarrow \LH(P \Box_G \E) $. \\
Now let $ (\E, \phi, F) $ be a $ G $-equivariant Kasparov $ A $-$ B $-module with $ G $-invariant 
operator $ F $. Since $ G $ is compact we may restrict to such Kasparov modules 
in the definition of $ KK^G $. Using our previous observations it follows easily that
$ (P \Box_G \E, \id \Box_G \phi, \id \Box_G F) $ is an $ H $-equivariant Kasparov $ P \Box_G A $-$ P \Box_G B $-module. 
It is not difficult to check that this assignment is compatible with homotopies and Kasparov products. \\
As a consequence we obtain the desired functor $ KK^G \rightarrow KK^H $
extending the functor $ F: G \Alg \rightarrow H \Alg $ defined above. 
This functor, again denoted by $ F $, preserves exact triangles and suspensions. \\
We may summarize our considerations as follows. 
\begin{theorem} \label{MEKK}
If $ H $ and $ G $ are monoidally equivalent compact quantum groups then 
the triangulated categories $ KK^H $ and $ KK^G $ are equivalent. 
\end{theorem}
Note that a trivial $ G $-$ C^* $-algebra is mapped to the corresponding trivial $ H $-$ C^* $-algebra
under this equivalence. As a consequence we immediately obtain the following assertion. 
\begin{theorem} \label{BCmon}
Let $ G $ and $ H $ be torsion-free discrete quantum groups whose dual compact quantum groups 
are monoidally equivalent. Then $ G $ satisfies the strong Baum-Connes conjecture iff $ H $ 
satisfies the strong Baum-Connes conjecture. 
\end{theorem}
\proof Using Baaj-Skandalis duality and theorem \ref{MEKK} we see that $ KK^G $ and $ KK^H $ are 
equivalent triangulated categories. In addition, the compactly induced $ G $-$ C^* $-algebra $ C_0(G) \otimes A $ corresponds to 
the compactly induced $ H $-$ C^* $-algebra $ C_0(H) \otimes A $ under this equivalence. Hence $ \bra \CI_G \ket = KK^G $ holds 
iff $ \bra \CI_H \ket = KK^H $ holds. \qed \\
Combining theorem \ref{BCmon} with theorem \ref{BCsuq2} yields the 
main result of this paper. 
\begin{theorem} \label{BCfree}
Let $ n > 2 $ and $ Q \in GL_n(\mathbb{C}) $ such that $ Q \overline{Q} = \pm 1 $. Then the free orthogonal quantum group $ \mathbb{F}O(Q) $ 
satisfies the strong Baum-Connes conjecture. 
\end{theorem}

\section{Applications} \label{secapp}

In this section we discuss consequences and applications of theorem \ref{BCfree}. In particular, we 
show that free orthogonal quantum groups are $ K $-amenable and compute their $ K $-theory. \\
The concept of $ K $-amenability, introduced by Cuntz for discrete groups in \cite{Cuntzkam},  
extends to the setting of quantum groups in a natural way \cite{VergniouxKamen}. 
More precisely, a discrete quantum group $ G $ is called $ K $-amenable if 
the unit element in $ KK^G(\mathbb{C}, \mathbb{C}) $ can be represented by a Kasparov module  
$ (\E, \pi, F) $ such that the representation of $ G $ on the Hilbert space $ \E $ is weakly contained 
in the regular representation. 
As in the case of discrete groups, this is equivalent to saying that 
the canonical map $ G \ltimes_\max A \rightarrow G \ltimes_\red A $
is an isomorphism in $ KK $ for every $ G $-$ C^* $-algebra $ A $. \\
Of course, every amenable discrete quantum group is $ K $-amenable. It is known \cite{Banicaunitary} that $ \mathbb{F}O(Q) $ is not amenable for 
$ Q \in GL_n(\mathbb{C}) $ with $ n > 2 $. \\
The main application of theorem \ref{BCfree} is the following result. 
\begin{theorem} \label{FOtheorem}
Let $ n > 2 $ and $ Q \in GL_n(\mathbb{C}) $ such that $ Q \overline{Q} = \pm 1 $. Then the free orthogonal quantum group $ \mathbb{F}O(Q) $ is 
$ K $-amenable. In particular, the map 
$$
K_*(C^*_\max(\mathbb{F}O(Q))) \rightarrow K_*(C^*_\red(\mathbb{F}O(Q))) 
$$
is an isomorphism. \\
The $ K $-theory of $ \mathbb{F}O(Q) $ is
$$
K_0(C^*_\max(\mathbb{F}O(Q))) = \mathbb{Z}, \qquad K_1(C^*_\max(\mathbb{F}O(Q))) = \mathbb{Z}.
$$
These groups are generated by the class of $ 1 $ in the even case and the class of the fundamental matrix $ u $ in the odd case. 
\end{theorem}
\proof Let us write $ G = \mathbb{F}O(Q) $. The reduced and full crossed product functors $ KK^G \rightarrow KK $ agree 
on $ \bra \CI \ket $ because they agree for all compactly induced $ G $-$ C^* $-algebras. Indeed, for $ B \in KK $ we have 
$$ 
G \ltimes_\max \ind_E^G(B) = G \ltimes_\max (C_0(G) \otimes B) \cong \KH \otimes B \cong G \ltimes_\red (C_0(G) \otimes B) 
= G \ltimes_\red \ind_E^G(B) 
$$ 
by strong regularity. According to theorem \ref{BCfree} it follows that the canonical map 
$ G \ltimes_\max A \rightarrow G \ltimes_\red A $ is an isomorphism in $ KK $ for every $ G $-$ C^* $-algebra $ A $, 
and this means precisely that $ G $ is $ K $-amenable. \\
As in section \ref{secbc} we denote by $ \mathfrak{J} $ the homological ideal in $ KK^G $ given by the kernel of the restriction 
functor $ \res^G_E: KK^G \rightarrow KK $. 
To compute the $ K $-groups $ K_*(C^*_\max(G)) \cong K_*(C^*_\red(G)) $ we shall construct a concrete $ \mathfrak{J} $-projective 
resolution of the trivial $ G $-$ C^* $-algebra $ \mathbb{C} $. \\
Let $ \hat{G} $ be the dual compact quantum group of $ G $, and let us identify the set of irreducible representations of 
$ \hat{G} $ with $ \mathbb{N} $. We write $ \pi_k $ for the 
representation corresponding to $ k \in \mathbb{N} $ and denote by $ \H_k $ the underlying Hilbert space. In particular, 
$ \pi_0 = \epsilon $ is the trivial representation. Moreover, $ \pi_1 $ identifies with
the fundamental representation given by $ u $, and we have
$ \dim(\H_1) = n $. The representation ring $ R(\hat{G}) $ 
of $ \hat{G} $ is isomorphic to the polynomial ring $ \mathbb{Z}[t] $ such that $ t $ corresponds to $ \H_1 $. \\
Due to the Green-Julg theorem and the Takesaki-Takai duality theorem, 
we have a natural isomorphism
$$
KK^{\hat{G}}(\mathbb{C}, G \ltimes_\red B) \cong K(\KH_G \otimes B) \cong K(B) 
$$
for $ B \in KK^G $. Consequently, taking into account $ KK^{\hat{G}}(\mathbb{C}, \mathbb{C}) \cong R(\hat{G}) $, 
the Kasparov product 
$$ 
KK^{\hat{G}}(\mathbb{C}, \mathbb{C}) \times KK^{\hat{G}}(\mathbb{C}, G \ltimes_\red B) 
\rightarrow KK^{\hat{G}}(\mathbb{C}, G \ltimes_\red B)
$$ 
induces an $ R(\hat{G}) $-module structure on $ K(B) $, and every element $ f \in KK^G(B, C) $ defines an 
$ R(\hat{G}) $-module homomorphism $ f_*: K(B) \rightarrow K(C) $. \\
For $ B = C_0(G) $ this construction leads to the action of $ R(\hat{G}) $ on itself by multiplication, and 
for $ B = \mathbb{C} $ the corresponding module structure on $ \mathbb{Z} $ is induced by the augmentation homomorphism 
$ \epsilon: \mathbb{Z}[t] \rightarrow \mathbb{Z} $ given by $ \epsilon(t) = n $. \\
Let us now consider the Koszul complex 
$$
\xymatrix{
0 \ar@{->}[r] & C_0(G) \ar@{->}[r]^{n - T} & C_0(G) \ar@{->}[r]^{\;\lambda} & \mathbb{C} 
}
$$
in $ KK^G $ defined as follows. The map $ \lambda: C_0(G) \rightarrow \KH_G \cong \mathbb{C} $ 
is given by the regular representation. Moreover, $ n: C_0(G) \rightarrow C_0(G) $ denotes the sum of $ n $ copies of the 
identity element. The morphism $ T: C_0(G) \rightarrow C_0(G) \otimes M_n(\mathbb{C}) \cong C_0(G) $ 
is the $ * $-homomorphism induced by the comultipliation 
$ \Delta: C_0(G) \rightarrow M(C_0(G) \otimes C_0(G)) $ followed by projection onto the matrix block corresponding to 
$ \pi_1 $ in the second factor. \\
Let us determine the map $ T_*: R(\hat{G}) \rightarrow R(\hat{G}) $ induced by $ T $ on the level of $ K $-theory. 
Consider the Hopf $ * $-algebra $ \mathbb{C}[\hat{G}] $ of matrix elements for $ \hat{G} $, and denote by 
$ \bra \;\,, \;\ket $ the natural bilinear pairing between $ \mathbb{C}[\hat{G}] $ and $ C_0(G) $. 
Under this pairing, the character $ \chi_k \in \mathbb{C}[\hat{G}] $ of the representation $ \pi_k $ corresponds to the trace on 
$ \KH(\H_k) $. Moreover we observe
$$
\bra \Delta(f), \chi_l \otimes \chi_1 \ket = \bra f, \chi_l \chi_1 \ket = \bra f, \chi_{l + 1} + \chi_{l - 1} \ket 
= \bra f, \chi_{l + 1} \ket + \bra f, \chi_{l - 1} \ket 
$$
for $ f \in \KH(\H_k) \subset C_0(G) $. According to the definition of $ T $ this implies
$$
T_*(\H_k) = \H_{k + 1} + \H_{k - 1}
$$
for all $ k \in \mathbb{N} $, where we interpret $ \H_{-1} = 0 $. It follows that 
$ T_* $ corresponds to multiplication by $ t $ under the identification $ \mathbb{Z}[t] = R(\hat{G}) $. \\
Applying $ K $-theory to the Koszul complex yields the exact sequence of $ \mathbb{Z}[t] $-modules
$$
\xymatrix{
0 \ar@{->}[r] & \mathbb{Z}[t] \ar@{->}[r]^{n - t} & \mathbb{Z}[t] \ar@{->}[r]^{\epsilon} & \mathbb{Z} \ar@{->}[r] & 0 
}
$$
where $ \epsilon $ is again the augmentation homomorphism given by $ \epsilon(t) = n $. Taking into account that 
$ C_0(G) $ is $ \mathfrak{J} $-projective, it follows easily that the Koszul complex yields a $ \mathfrak{J} $-projective 
resolution of $ \mathbb{C} $. \\
Since the Koszul resolution has length $ 1 $, we obtain a Dirac morphism $ \tilde{\mathbb{C}} \rightarrow \mathbb{C} $ 
as in the proof of theorem 4.4 in \cite{MNhomalg1}. The only piece of information that we need about this construction is that the 
$ G $-$ C^* $-algebra $ \tilde{\mathbb{C}} $ fits into an exact triangle 
$$
\xymatrix{
C_0(G) \ar@{->}[r] & C_0(G) \ar@{->}[r] & \tilde{\mathbb{C}} \ar@{->}[r] & \Sigma C_0(G)
}
$$
in $ KK^G $. Here the first arrow $ C_0(G) \rightarrow C_0(G) $ is given by $ n - T $, but 
we will not make use of this fact. 
By applying the crossed product functor we obtain an exact triangle 
$$
\xymatrix{
\KH \ar@{->}[r] & \KH \ar@{->}[r] & G \ltimes_\max \tilde{\mathbb{C}} \ar@{->}[r] & \Sigma \KH
}
$$
in $ KK $. Hence the associated long exact sequence in $ K $-theory takes the form 
$$
\xymatrix{
 {\mathbb{Z}\;} \ar@{->}[r] \ar@{<-}[d] &
      K_0(G \ltimes_\max \tilde{\mathbb{C}}) \ar@{->}[r] &
        0 \ar@{->}[d] \\
   {\mathbb{Z}\;} \ar@{<-}[r] &
    {K_1(G \ltimes_\max \tilde{\mathbb{C}})}  \ar@{<-}[r] &
     {0} \\
}
$$
Since $ C^*_\max(G) $ has a counit the group
$ K_0(G \ltimes_\max \tilde{\mathbb{C}}) \cong K_0(C^*_\max(G)) $ contains a direct summand $ \mathbb{Z} $ generated 
by the unit element $ 1 \in C^*_\max(G) $. It follows that the upper left horizontal map in the diagram is 
an isomorphism. Hence the vertical arrow is zero, and the lower left horizontal map is an isomorphism as well. \\
It remains to identify the generator of $ K_1(C^*_\max(G)) \cong K_1(G \ltimes_\max \tilde{\mathbb{C}}) $. 
Clearly, the fundamental unitary $ u \in M_n(C^*_\max(G)) $ defines an element $ [u] \in K_1(C^*_\max(G)) $. 
The discussion at the end of section 5 in \cite{BdRV} shows that we find a quotient homomorphism 
$ \pi: C^*_\max(G) \rightarrow C^*_\max(\mathbb{F}O(M)) $ 
for some matrix $ M \in GL_2(\mathbb{C}) $. On the level of $ K $-theory, the class 
$ [u] $ maps under $ \pi $ to the class of the fundamental matrix of $ C^*_\max(\mathbb{F}O(M)) $ 
in $ K_1(C^*_\max(\mathbb{F}O(M))) $. Since $ M \in GL_2(\mathbb{C}) $, 
the quantum group $ \mathbb{F}O(M) $ is isomorphic to the dual of $ SU_q(2) $ for some $ q \in [-1, 1] \setminus \{0\} $. \\
Let $ u_q \in M_2(C(SU_q(2))) $ be the fundamental matrix of $ SU_q(2) $. For positive $ q $,  
the index pairing of $ u_q $ with the Fredholm module corresponding to the Dirac operator
on $ SU_q(2) $ is known to be equal to $ 1 $, see \cite{DLSVV}. The argument extends to the case of 
negative $ q $ in a straightforward way. \\
As a consequence of these observations we conclude that $ [u] $ is a generator of $ K_1(C^*_\max(G)) = \mathbb{Z} $, 
and this finishes the proof. \qed \\ 
In his work on quantum Cayley trees \cite{Vergniouxtrees}, Vergnioux has constructed an analogue of the Julg-Valette element for 
$ \mathbb{F}O(Q) $. Our considerations imply that this element is homotopic to the 
identity, although we do not get an explicit homotopy. 
\begin{cor} 
Let $ n > 2 $ and $ Q \in GL_n(\mathbb{C}) $ such that $ Q \overline{Q} = \pm 1 $. 
The Julg-Valette element for $ \mathbb{F}O(Q) $ is equal to $ 1 $ in 
$ KK^{\mathbb{F}O(Q)}(\mathbb{C}, \mathbb{C}) $. 
\end{cor}
\proof Let us write $ G = \mathbb{F}O(Q) $. 
An analogous argument with $ K $-homology instead of $ K $-theory as in the proof of theorem \ref{FOtheorem} 
shows that $ KK(C^*_\max(G), \mathbb{C}) \cong \mathbb{Z} $ is generated by the class of the counit 
$ \epsilon: C^*_\max(G) \rightarrow \mathbb{C} $. It follows that the forgetful map 
$ KK^G(\mathbb{C}, \mathbb{C}) \rightarrow KK(\mathbb{C}, \mathbb{C}) \cong \mathbb{Z} $ 
is an isomorphism. Since the Julg-Valette element in \cite{Vergniouxtrees} has index $ 1 $ this yields the claim. \qed \\
In the special case $ Q = 1 \in GL_n(\mathbb{C}) $ the quantum group $ \mathbb{F}O(Q) = \mathbb{F}O(n) $ is unimodular. 
Hence the Haar functional $ \phi: C^*_\red(\mathbb{F}O(n)) \rightarrow \mathbb{C} $ is a trace and determines an additive map 
$ \phi_0: K_0(C^*_\red(\mathbb{F}O(n))) \rightarrow \mathbb{Z} $. 
\begin{theorem}
Let $ n > 2 $. The free quantum group $ \mathbb{F}O(n) $ satisfies the analogue of the Kadison-Kaplansky conjecture. 
That is, $ C^*_\red(\mathbb{F}O(n)) $ does not contain nontrivial idempotents. 
\end{theorem}
\proof The classical argument for free groups carries over. Since every idempotent is similar to a projection it suffices to show 
that $ C^*_\red(\mathbb{F}O(n)) $ does not contain nontrivial projections. 
We know that the Haar functional is a faithful tracial state on 
$ C^*_\red(\mathbb{F}O(n)) $. 
Assume that $ p \in C^*_\red(\mathbb{F}O(n)) $ is a projection. Then from the positivity of $ \phi $ 
we obtain $ \phi(p) \in [0,1] $, and from the above 
considerations we know $ \phi(p) = \phi_0([p]) \in \mathbb{Z} $. This implies $ p = 0 $ or $ 1 - p = 0 $. \qed \\
Finally, we note that the dual of $ SU_q(2) $ does not satisfy the analogue of the Kadison-Kaplansky conjecture. 
In fact, there are lots of nontrivial idempotents in $ C(SU_q(2)) $ for $ q \in (-1,1) \setminus \{0\} $. 

\bibliographystyle{plain}

\end{document}